\documentclass{article}

\usepackage{float}
\usepackage{amsfonts, enumerate} 
\usepackage{graphicx, graphics,stmaryrd} 

\usepackage{lipsum}

\usepackage{epstopdf}
\usepackage{algorithmic}
\ifpdf
  \DeclareGraphicsExtensions{.eps,.pdf,.png,.jpg}
\else
  \DeclareGraphicsExtensions{.eps}
\fi

\usepackage{amscd}

\usepackage{tensor}
\usepackage[english]{babel}
\usepackage[utf8]{inputenc}
\usepackage{graphics}
\usepackage{amsmath}
\usepackage{amstext}
\usepackage{engrec}
\usepackage{rotating}
\usepackage[safe,extra]{tipa}
\usepackage{multirow}
\usepackage{enumitem}
\usepackage{bm}
\usepackage{stackrel} 
\usepackage{tikz-cd} 
\numberwithin{equation}{section}

\newcommand{\RP}{P}

\def\R{\mathbb{R}}
\def\Z{\mathbb Z}
\DeclareMathOperator{\id}{id}

\newcommand{\mc}{\mathcal}

\newcommand{\m}{\mbox}
\newcommand{\norm}[1]{\left\lVert#1\right\rVert}

\newtheorem{teo}{Theorem}[section]
\newtheorem{prop}[teo]{Proposition}
\newtheorem{defin}[teo]{Definition}
\newtheorem{rmk}[teo]{Remark}
\newtheorem{ex}[teo]{Example}

\title{The cortical V1 transform as a heterogeneous Poisson problem}
\author{A.Sarti, M.Galeotti, G.Citti} 

\usepackage{amsopn}

\begin{document}

\maketitle


\begin{abstract}

Receptive profiles of the primary visual cortex (V1) cortical cells are very heterogeneous and act by differentiating the stimulus image as operators changing from point to point. In this paper we aim to show that the distribution of cells in V1, although not complete to reconstruct the original image, is sufficient to reconstruct the perceived image with subjective constancy. We show that a color constancy image can be reconstructed as the solution of the associated inverse problem, that is a Poisson equation with heterogeneous differential operators.  At the neural level the weights of short range connectivity constitute the fundamental solution of the Poisson problem adapted point by point. A first demonstration of convergence of the result towards homogeneous reconstructions is proposed by means of homogenization techniques. 

\end{abstract}



\section{Introduction}

The visual brain extracts more and more complex features starting from the visual stimulus. The cells of the retina and the lateral geniculate nucleus (LGN) extract the position of the contours, then the cells of the primary visual cortex (V1) extract the position and orientation of the image contours, up to the upper cortices where more complex  features are extracted.
The action of cells is characterized to a first approximation by linear receptive profiles (Ps) \cite{HW59}, called also ``classical receptive profiles", which behave like filters
representing the impulse response of the cells.
Predictions of receptive profiles in V1 are in good agreement with receptive measurements reported in the literature (Hubel and Wiesel  \cite{HW59,HW62,HW05}; 
DeAngelis et al.~\cite{DA95,DA04}). Specifically, explicit  neurophysiological caracterisations have been given of LGN neurons in terms of Laplacian of Gaussian,
 and simple cells in V1 have been compared to related models in terms of Gabor functions (Marcelja \cite{MA80}; 
Jones and Palmer\cite{JP87,JP872}), differences of Gaussians (Rodieck \cite{R65}) or Gaussian derivatives (Koenderink 
and van Doorn \cite{K87}; Young  \cite{Y87}; Young et al. \cite{YLM01,YL01}; Lindeberg \cite{TL13}). 

In several papers it has been shown that the set of profiles can be obtained through group transformations starting from a mother profile. For example in \cite{CS06} simple cell profiles of V1 were studied in the rotation and translation group and it was shown that all the functional architecture related to these cells (horizontal connectivity and association fields) is to be considered as a direct consequence of the Lie group symmetries involved. This finding was extended to cells sensitive to scale  in the symplectic group \cite{SCP08}, to movement in the Galilean group \cite{BA13}, to curvature in the Engel group \cite{A12} and to rotation, scale and frequency in the Heisenberg group  \cite{BSC20}. The interested reader could refer to  \cite{CS14} for a complete review.

Simple cells sensitive to position and orientation are topographically organized in V1. The process of visual mapping from the retina to cortical neurons is known as retinotopy. Moreover in primates and cats, neurons with similar orientation selectivity are clustered together to constitute the so-called hypercolumnar \cite{HW59} or pinwheel  \cite{BG91} organisation. In the opposite, rodents do not have such specific feature organisation as primates and cats. This means that the information encoded in V1 in rodents is distributed in a less organised fashion, known as ‘salt-and-pepper maps’   \cite{MK14}. 
Anyway,  whatever the feature organization is, the 3-dimensional group of rotations and translations of the plane is not completely represented since it is projected in the 2D cortical surface of V1. 

An important problem of contemporary neuroscience consists in understanding whether the perceived image can be reconstructed  starting from the partial information carried by cells in V1. In other words, it is a question to understand whether V1 is a high-resolution image buffer that contains the  representation of the perceived image. In this sense V1 would make a cortical transform of the stimulus to obtain the perceived image.

In this regard, Davide Barbieri in \cite{DB21} shows that the distribution of receptive profiles in V1 is not complete, therefore it is not sufficient to reconstruct the visual stimulus. Subsequently, he shows which additional constraints would be required to achieve the reconstruction.

In this paper we aim to show that the distribution of cells in V1, although not complete to reconstruct the original image, is sufficient to reconstruct the perceived image without additional constraints. In particular we focus on the perceptual phenomena of lightness constancy and color constancy, which are fundamental characteristics of the perceived image.   Lightness constancy refers to the perception of an object’s lightness as invariant with respect to the illumination conditions.
Analogously, color constancy is a similar phenomenon but this time is the perception of color resulting invariant with respect to illumination.
We therefore aim to show that the distribution of cells in V1 is sufficient to perform a transform of the visual stimulus in the   lightness and color constancy image perceived by the subject.


To do this we will consider the Ps as Gaussian derivatives \cite{K87, Y87, YLM01,YL01}  with heterogeneous metrics and order of derivation. They act by differentiating the image as heterogeneous operators changing from point to point. Then we show that an image can be reconstructed as the solution of the associated inverse problem, that is a Poisson equation with heterogeneous operators.  The homogeneous equivalent of this reconstruction is the Retinex differential model presented in \cite{Ki03, Mo11} in which the perceived image was obtained by first differentiating the stimulus with homogeneous Laplacian operators and then solving the associated inverse problem (Poisson equation). In our case we will generalize the Poisson equation to equations changing from point to point in order to reconstruct the perceived image.  This kind of differential problems with heterogeneous operators was introduced, among others, in \cite{SCP19, SCP22} to define the possibility of performing a differential calculus with operators that change from point to point in space and time.

Observe that the heterogenous differential structure of images has been largely studied in the past, starting from the pioneristic work of Jan Koenderink \cite{K84}. A wide amount of studies about locally adaptive frames (or Gauge frames) has been developed and an exhaustive list of references has been collected in  \cite{HR03}. 
In this setting, detectors of features like edges, corners, t-junctions, ridges, monkey-saddles and many more,e have been modeled by heterogeneous differential operators. A further expansion of heterogeneity has been provided by data-driven left-invariant metrics \cite{Be15}.  These works, although they underline the heterogeneity of the differential structure of  images, concern the direct problem of differentiation and feature extraction rather than the inverse problem of stimulus reconstruction.  On the other hand, 
research about stimulus reconstruction by inversion of differential operators are very few, and up to now they deal with homogeneous operators. \cite{M-G18, BGS,BMSC}. 

The originality of the present research is in performing operator inversion in a differential heterogeneous setting.
We provide an existence results for weak solutions of the operator. 
A first proof of convergence of the result towards a homogeneous reconstruction as in \cite{Ki03, Mo11}, will be proposed by means of homogeneisation techniques.
The proof is constructed with operators discretised on regular grids (difference operators). 
We refer in particular to 
\cite{KKG82, koz87, piatremy01}, and
 develop the notion of H-convergence (see Definition \ref{410}) that was initially 
introduced by Spagnolo \cite{Spagnolo} and De Giorgi \cite{De Giorgi75, De Giorgi83}. Results of this kind for mathematical models of the cortex are for example \cite{SCM03,CMS03}.
We consider a family of difference operators with random coefficients on an $\varepsilon$-spaced
reticulum, and we prove that it can converge 
as $\varepsilon \to 0$ to a deterministic differential
operator  defined on a continuous domain
and isotropic under suitable hypothesis. 

The paper is organized as follows. In Section \ref{RP} some classical concepts of functional architecture of V1 are recalled. 
Classical models of receptive profiles are presented in terms of Gaussian derivatives and  orientation maps of such profiles 
are shown both for hypercolumnar organisation and salt-and-pepper maps. 
In Section \ref{CT}, the main model of cortical transform as a heterogeneous Poisson problem is presented.  
In Section \ref{existence} we prove the existence of a weak solution of the problem, using a steepest descent method. 
In the subsequent Section  \ref{DC}, the notion of $H$-convergence is introduced and a sketch 
of the proof of convergence of the heterogeneous problem towards the homogeneous one is proposed. 
The main numerical results are shown and discussed in Section \ref{NR}.

\section{Functional architecture of V1} \label{RP}

\subsection{Receptives profiles}

A receptive profile is the impulse response of a cell, that is the response of the cell to a delta of Dirac of retinal stimulus.  Statistical studies on the Ps of V1 cells in macaques  show a great heterogeneity of behaviours  \cite{Ringach}.  Copies of  center-surrounds Ps (Mexican hats) are present, as well as a great variety of simple cells with anisotropic Ps that detect orientation of boundaries, all with different orientation and frequency and many others complex cells.

In the simple case of Mexican hats the receptive profile 
will be denoted $P_h$ (the index $h$ denotes the fact that $P$ is the profile of a Mexican hat), and it
is axial symmetric and in terms of gaussian derivatives writes:

\begin{equation}
\label{mh}
\RP_h(x_1, x_2) =\Delta  G_\sigma(x_1, x_2),
\end{equation}
where $x=(x_1, x_2)$ is the general point of the plane
where $G_\sigma$ is a two-dimensional Gauss function with variance $\sigma$ and $\Delta$ is the Euclidean Laplacian. 

\begin{figure}[htbp]
\centering
\label{f1x}
\includegraphics[height= 1.8in ]{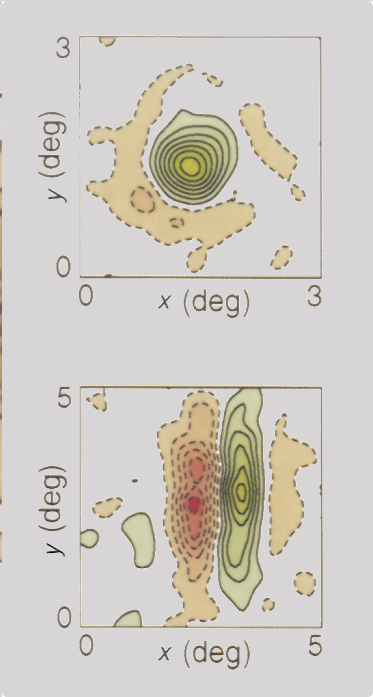}
\includegraphics[height= 1.8in ]{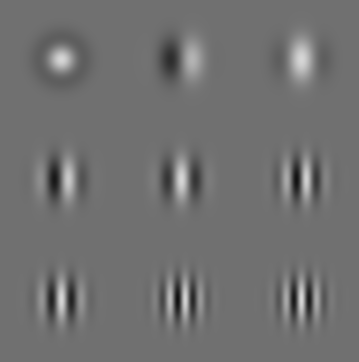}

\caption{Left: Receptive profiles of different kind of cells after neurophysiological measurement (inhibitory response in red, excitatory in yellow) as reported in  \cite{DA95}\cite{DA04}. Upper Left: LGN receptive profiles show a characteristic center-surround response that can be modelled as a Laplacian of Gaussian. Bottom Left: Simple cells of V1 show strong orientation preference and are modelled by directional derivatives of Gaussian. Right: Representation of Ps as derivatives of Gaussian of different degrees.}
 \label{f1} 
\end{figure}

In case of receptive profiles of simple cells let's consider the directional derivative  
\begin{equation}\label{defX1}
X_{\theta, 1}=\cos \theta \partial_{x_1}+\sin\theta \partial_{x_2}.
\end{equation}  The receptive profile of a simple cell with preferred orientation $ \theta+\pi/2$ writes:

\begin{equation}
\label{sc}
\RP_s(x_1, x_2)={{X^ {2\beta}_{\theta, 1}}} G_\sigma(x_1, x_2), 
\end{equation}
The coefficient $\beta$ is always an integer,
so that ${{X^{2\beta}_{\theta,1}}}$ denotes the directional derivative of order $2\beta$ in the direction $\theta$. In particular this is a derivative of even order. Tipically we will be interested in $\beta=1$ or $\beta=2$, giving rise to derivatives of order $2$ or $4$. 
The higher order derivative denotes the presence of a higher number of sinusoidal cycles under the Gauss bell 
\cite{K87} (Figure \ref{f1}). 

\subsection{Orientation maps}
\begin{figure} [htbp]
\label{f2x}
\centering
\includegraphics[height= 1.6in ]{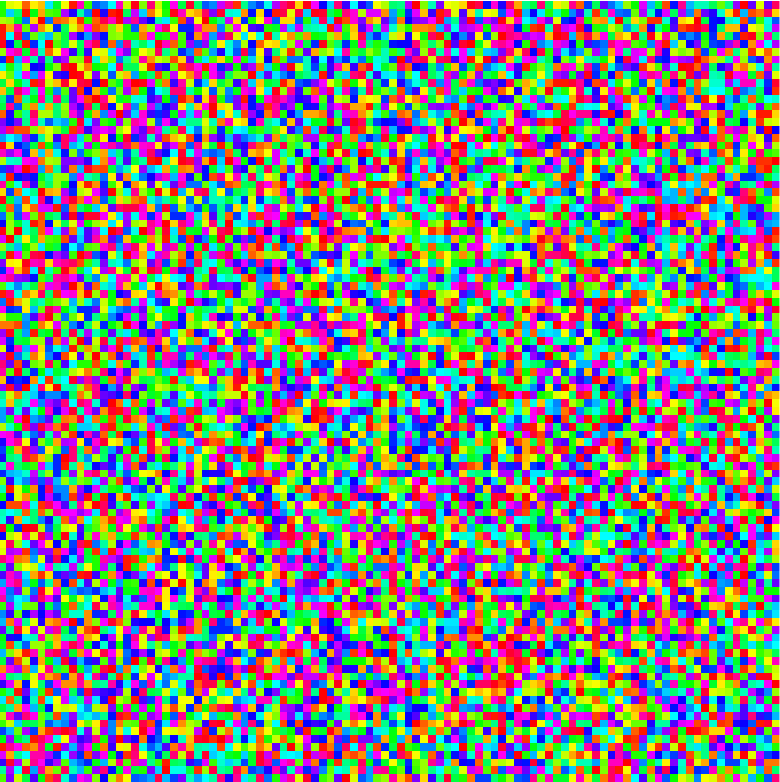}
\includegraphics[height= 1.6in ]{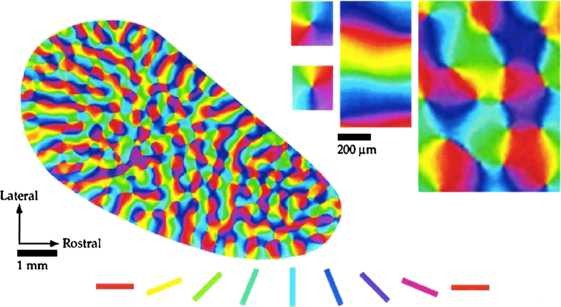}
\caption{Orientation maps in rodents are salt and pepper noise (left, our reconstruction) while in primates they present the so-called pinwheels structure  (right, reproduced from \cite{BO97}). A mathematical model of the two distributions is given respectively by Equation \eqref{salt}  and Equation~\eqref{pin}. }
\label{f2}
\end{figure}

The preferred orientation of simple cells changes in the cortex point to point  \cite{HW59} giving rise to the so called orientation map. The orientation map is a map $\theta:\mathbb{R}^2\rightarrow [0,\pi]$, where $x =(x_1, x_2)$ are cortical coordinates and $\theta= \theta(x_1,x_2)$ is the preferred direction of columns of simple cells. In primates and cats, neurons with similar orientation selectivity are clustered together to constitute the so-called pinwheel  organisation  \cite{BG91}. 
A simple model of pinwheel shaped orientation map is proposed in \cite{Petitot} where the
map is obtained through the superposition of randomly
weighted complex sinusoids:

\begin{equation}
\label{pin}
\theta_{1}(x_1, x_2)=\arg \sum_{k=1}^N c_k e^{i2\pi(x_1\cos(2\pi k/N)+x_2\sin(2\pi k/N))}
\end{equation}
with $N$ denoting the number of orientation samples and
where the coefficients $c_k \in [0, 1]$ are  white noise. From equation \ref{pin} it is obtained the so called pinwheel structure, that is an orientation map with the presence of singularities around which all the orientations are present.  The function is complex valued and the color maps the argument of complex variable. The distance between singularities is almost constant. Pinwheel points are the two dimensional implementation of a hypercolumn of orientations that is the set of orientations corresponding to a retinal point in the original model of Hubel and Wiesel \cite{HW59} (Figure \ref{f2}).

In the opposite in rodents the orientation information encoded in V1 is distributed in a less organised fashion, known as ‘salt-and-pepper maps’   \cite{MK14} and can be modelled by: 
\begin{equation}
\label{salt}
\theta_{2}(x_1,x_2)= d(x_1, x_2)
\end{equation}
where  the coefficients $d \in [0, \pi]$ are white noise.

\subsection{Short and long range connectivity}
\begin{figure}[htbp]
\label{f3x}
\centering
\includegraphics[height= 1.5in ]{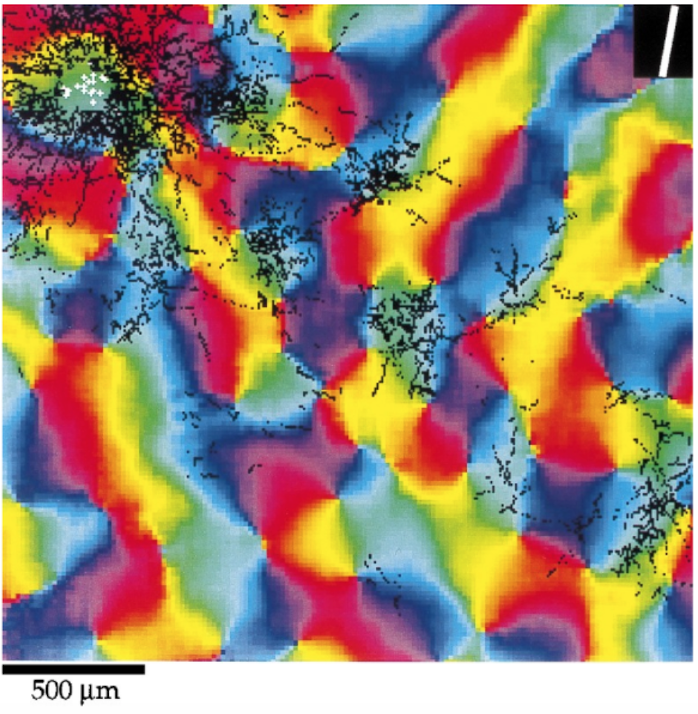}
\includegraphics[height= 1.5in ]{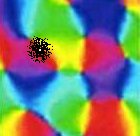}

\caption{Left: Long range connectivity between hypercolumns from \cite{BO97}. Right: Short range connectivity whithin a hypercolumn (our representation from \cite{DG}). } 
\label{f3}

\end{figure}

Cortical cells are reciprocally connected by short and long range connectivity. Short range connectivity corresponds to the distribution of interaction between cells within a hypercolumn. On the other hand long range connectivity induces interactions between hypercolumns. Long range connectivity is very anisotropic and connects mainly cells with a similar orientation preference, while short range connections are isotropic.  
Particularly, from neurophysiological measurements of the activity of pairs of V1 neurons it is computed a correlation function that is proportional to the functional connectivity between the pairs of units. This short range connectivity show radial symmetry and a decreasing in intensity proportional to  
$\log(r )$ or $1/r$   \cite{DG}  (Figure \ref{f3}).

These connections are modulatory, meaning that they act on the output of Ps of cortical cells and not directly on the  thalamic stimulus.

\section{ The cortical transform} \label{CT}
%

We describe in the following the model for the coupled activity of classical receptive profiles and the short range connectivity that modulates the feed forward input giving rise to contextual modulation effects. Particularly the weights of local connectivity constitute the fundamental solution of a Poisson problem where the forcing term is given by the output of classical receptive profiles.   This coupling can account for the lightness constancy of the perceived image (Retinex effect). Given the heterogeneity of cells linked by local connectivity, the model will end up in a Poisson equation with heterogeneous operators. Homogeneous version of the model can be found for example in \cite{Ki03}  \cite{Mo11}   while a more complex model but without any image reconstruction is proposed by  Bressloff and Cowan  \cite{bressloff} . 

\subsection{Feed forward action on the input stimulus}

The action of a classical receptive profile $\RP$ on an input stimulus $I$ writes:
$$
O=\RP \star I
$$
where  $\star$ is the convolution product and $I$ is the image stimulus, on which the P acts as a linear filter.
Particularly the action of Mexican reads
\begin{equation}\label{euclidean}
O_h=\RP_h \star I = \Delta  G_\sigma \star I = \Delta I_\sigma, 
\end{equation}
where $ I_\sigma =  G_\sigma \star I$ is smooth approximation of $I$. 
The action of the simple cells $\RP_s$ writes:
\begin{equation}\label{sub}
O_s=-(-1)^\beta X_{\theta,1}^{2\beta} G_\sigma \star I= -(-1)^\beta X_{\theta, 1}^{2\beta} I_\sigma,
\end{equation}
where $\beta$ is a integer positive and even.
For  $\beta=1$ the operator ${X}_{\theta, 1}^{2\beta}$ is a second order directional derivative in the direction $\theta,$ which can be identified with a degenerate Laplacian. Let us recall that $-\Delta $ and $-{X}_{\theta, 1}^{2}$ are elliptic.

Hence we postulate that the operator which acts on the image is different from a point to the other 
and can be a random combination of the previous operators. 
In the sequel we will always denote by $L_\Lambda$ the considered linear operator which can be 
$
L_\Lambda = \Delta \; \text{or }L_\Lambda= {X}_{\theta, 1}^{2\beta}, 
$
or a linear combination of them. We will be 
mostly interested in the values $\beta=1$ and $\beta=2$, which denote derivatives of order 2 or 4 respectively. The general operator will  read
\begin{equation}\label{operats}
L_\Lambda = a_1 \Delta + a_2{X}_{\theta, 1}^{2} - a_3   (-1)^\beta{X}_{\theta, 1}^{2\beta}
\end{equation}
for a partition of the unity $a_i $: at every point $x$ one of these coefficients will take the value 1, and the others will vanish.

\subsection{Reconstruction of the perceived image} 
The output of a distribution of V1 cells in response to the visual input $I$ is the differentiation of the visual input itself. We will propose here a  procedure to construct the perceived image by considering the integration of the operator.
In previous papers, we conjectured that this step is accomplished by the action of the connectivity between V1 and LGN in a homogeneous setting  \cite{BGS,BMSC}. Here we clarify that it can be accomplished by short range connectivity of V1 and by considering heterogeneous operators of V1.

Let's consider the inverse problem with respect to the differentiation $ L_{\Lambda} I(x_1,x_2)$, that
is  the heterogeneous Poisson problem, with Neumann or Dirichlet conditions,

\begin{equation} \label{maineq}
L_{\Lambda}  u(x_1, x_2) =L_{\Lambda} I(x_1, x_2)
\end{equation}
where the right side term is known and represents the action of receptive profiles on the visual input, 
while the solution $u$ represents the reconstructed-perceived image. Equation (\ref{maineq}) can be considered as 
an extension  of a Poisson equation to heterogeneous operators. This process constructs the perceived image, 
that is the original input $I(x_1,x_2)$ up to a harmonic function. 
Clearly, harmonic means no more the annulation of the Laplacian but the annulation of the heterogeneous operator $L_{\Lambda}$.

Having formalized the expression of the operator $L_\Lambda$ we 
study the associated equation 
 (\ref{maineq}). Calling $f=  L_{\Lambda}I$ the equation reduces to 
 \begin{equation}\label{main2}
  L_{\Lambda} u = f
 \end{equation}
with Dirichlet or Neumann boundary conditions. 
We will prove the existence of solutions via parabolic approximation. Indeed we consider  
 the evolution equation 
 \begin{equation}\label{claprob1}
\partial_t \tilde u = L_{\Lambda}\tilde u - f
\end{equation}
with an arbitrary initial condition, and the same boundary condition as in (\ref{main2})  at any instant of time. 
Formally, for any open bounded domain $Q$
we will introduce suitable Sobolev spaces $W^{1, 2}_\Lambda(Q)$, we will consider $L_{\Lambda}$ as an operator acting on the space and we will look at  
the Cauchy problem 
\begin{align}\label{Cauchy}
&\tilde u' = L_{\Lambda}\tilde u - f\\
&u(0) \in W^{1, 2}_\Lambda(Q)\nonumber
\end{align}
(the initial condition vanishes at the boundary if we consider a Dirichlet problem). 
We will see under which condition the Cauchy problem is defined for every instant of time, there exists $\lim_{t \to \infty} \tilde u $ and the limit coincides with the solution $u$ of (\ref{main2}).

\section{Existence of weak solutions}\label{existence}
The operator defined in \eqref{operats} is degenerate
when $a_1=a_3=0$, even if $\theta$ is constant. 
Indeed, if $\theta=0$ the equation reduces to 
$$ \partial^2_{x_1x_1} u =f$$
in an open set $Q\subset \R^2,$ so that it is not uniformly elliptic.   
In addition it is of order 2 at some points, and of order 4 at other points. We will formalize the definition of solutions
 in suitable Sobolev spaces associated to the considered directional derivatives. 
In addition, we will prove the convergence of the solution of the parabolic 
equation \eqref{claprob1} to the solution to the problem \eqref{main2}.

\subsection{Definition of weak solution}
We would like to introduce a definition of a differential operator sufficiently flexible to allow to cover the Euclidean Laplacian and the degenerate Laplacian defined in \eqref{operats}. 
To this end we will consider a smooth manifold $Q$ locally diffeomorphic to $\R^2$, 
and we will define a sub-bundle 
$H$ of the tangent bundle $TM$ of $Q$. 
Then we define a metric $g$ on $H$. The triple $(Q, H, g)$ is called sub-Riemannian manifold. 
By definition of metric, for every $x$ the metric $g$ defines a scalar product on the horizontal tangent plane $H_{x}$.

\subsection{Second order operators}\label{secordop}

At every point $x$ the horizontal tangent space $H_x$ will be defined by the choice of generators $\Lambda_x$. 
The expression of $\Lambda_x$ will reflect the choice of receptive profiles at that point. 
As a consequence we will choose as set of generators 
 \begin{equation}\label{soliti2} \Lambda_{x} = \{z_1, z_2\},\text { where }z_1= a_1(\cos(\theta(x)), \sin(\theta(x)) ),\;  z_2 = a_2(-\sin(\theta(x)) ,\cos(\theta(x))). 
 \end{equation}
The dimension of the generated horizontal tangent space will depend on the zeros of the functions $a_i$, and can be 1 or 2. 
Correspondingly, in addition to the vector field $X_{\theta(x), 1}$ introduced in \eqref{defX1} we define a second vector field, to complete the basis of the space: 
\begin{equation}\label{vectors2}
X_{\theta(x), 1} =  \cos(\theta(x)) \partial_{x_1} + \sin(\theta(x))\partial_{x_2} ,\;  X_{\theta(x), 2} = -\sin(\theta(x)) \partial_{x_1} + \cos(\theta(x))\partial_{x_2}
\end{equation}
and we will work with the vector fields $a_i X_{\theta(x), i}.$ \newline

For any bounded set $Q\subset \R^d$, and every $h \in C_0^\infty$, we define the horizontal gradient 
$$\nabla_\Lambda h := \sum_{i=1}^2 X_{\theta(x), i}h X_{\theta(x), i}.$$
This allows to define the 
integral norm  
$$\norm{h}_{W^{1, 2}_\Lambda(Q)} = \sqrt{\sum_{i, j=1}^2 \int _Q  a_i^2({X_{\theta(x), i}} h)^2(x)  dx}.$$

The Sobolev space
$W^{1, 2}_{0, \Lambda}(Q)$ will be the closure of $C^\infty_0(Q)$ with respect to this norm.
If $Q$ is smooth,  $W^{1, 2}_{\Lambda}(Q)$ is the set of function $u$ whose  norm $\norm{h}_{W^{1, 2}_{\Lambda}(Q)} $ is bounded. Recall that the coefficients $a_i$ can vanish, making the space degenerate.
By definition the following Dirichlet functional is well defined in the space $W^{1, 2}_{\Lambda}(Q)$: 
\begin{equation}\label{defF}
F(u) = \sum_{i=1}^2 \int _Q  a_i^2({X_{\theta(x), i}} u)^2(x)  dx.
\end{equation}

If $X$ is an arbitrary operator of the second order with differentiable coefficients,
 we call formal adjoint of ${X}$
the operator $X^*$  defined by the relation
$$\int_Q {X} u h =  \int_Q  u X^* h \quad \  \forall u \in W^{1, 2}_{\Lambda}(Q), h\in C_{0}^\infty(Q).$$
The associated Laplacian is formally defined as 
\begin{equation} \label{op}
L_{\Lambda} u =\sum_{i =1}^2 X_{\theta(x), i}^* ( a^2_{i} { X_{\theta(x), i}} u).
\end{equation}

By definition of formal adjoint, the weak solution is defined as follows.
\begin{defin}
Let $Q$ be a bounded open set. 
A function $u\in W^{1, 2}_{\Lambda}(Q)$ is a weak solution of the equation $L_{\Lambda } u = f \in L^2(Q) $ 
if for every smooth function $h$ compactly supported in $Q$ 
$$ 
\int_{Q} \sum_{i =1} ^2a^2_{i}X_{\theta(x), i} u X_{\theta(x), i} h  dx =- \int_Q f h dx .$$ 
where the $a_i$ are bounded and measurables. 
The function $u$ attains the value $0$ at the boundary in a weak sense if $u\in W^{1, 2}_{0, \Lambda}(Q)$. It attains the Neumann boundary datum if it is the minimum of the associated functional 
\begin{equation}\label{J}
J(u) = F(u) - \int_Q f u.
\end{equation}

\end{defin}

\subsection{Higher order operators} 

We already noted that the operator defined in \eqref{operats} can have order grater than $2$ or can be a linear combination of operator of different order. 
In the simplest case  the operator reduces to 
\begin{equation}\label{general}
L_{ \Lambda}u(x) =  X_{\theta(x), 1}^*(a^2_1X_{\theta(x), 1} )u(x)  - ( X^2_{\theta(x), 1})^*(a^2_3 X^2_{\theta(x), 1} ) u(x).
\end{equation}
In order to simplify the proof, we will consider the operator
\begin{equation}\label{general}
L_{ \Lambda}u(x) =  X_{\theta(x), 1}^*(a^2_1X_{\theta(x), 1} u)(x)  - \Delta (a^2_3 \Delta) u(x)\end{equation}
that  can have different order at different points. 
Observe that $L_\Lambda=-dF$ where $F$ is the operator
\begin{equation}\label{Fu}F(u) = 
 \sum_{i=1}^2 \int _Q  a_i^2({X_{\theta(x), i}} u)^2(x)  dx +\int _Q  a^2_3(\Delta ^2 u)^2(x)  dx .
 \end{equation}
As before this functional, well defined on functions $u\in C^{\infty}_0$, defines the square of a norm. 
The closure of  $ C^{\infty}_0$  with respect to this norm, defines a Sobolev space where the associated equation is well defined.
If $a_3=0$, the operator reduces to a second order operator, which can be degenerate or not. 
If $a_{i}=a_2=0$, $a_3=1$, the operator reduces to the bi-Laplacian, denoted by $\Delta_E^2$.
We have defined an operator $L_{\Lambda}$ heterogeneous both in its differential order  and in the type of metric.

\subsection{Existence of a solution}

In all the considered problems, the operator $L_{\Lambda}$ is the differential of a functional $J$ defined on a suitable Hilbert  space $S$ of Sobolev type. 
We will show that in this case the solution of the stationary problem  can be obtained via the steepest descent method. 

Precisely assume that $Q$ is a bounded domain subset of $\R^2$. 
Also assume that $J$ is of class $C^{1,1}(S)$
By the Riesz representation theorem there exists a function 
$\nabla J(u) \in S$ such that $dJ(u)(h)$ admits the following representation 
with respect to the scalar product of the space $S$
$dJ(u)(h)=\langle \nabla J(u), h\rangle$. 
Assume that $Lu-f = - \nabla J(u)$, so that  the problem \eqref{Cauchy}
reduces to 
\begin{align}\label{Ca}
&\tilde u' =- \nabla J(\tilde u)\\
&\tilde u(0) \in S\nonumber
\end{align}

\begin{rmk}If $J\in C^{1,1}$ and it is bounded below, then
$- \nabla J(u)$ is Lipchitz continuous, and the Cauchy problem \ref{Ca} is well defined. Due to the boundness from below of $J$, its solution is defined on $[0, \infty[$. 
\end{rmk}
We will need the following compactness assumption of Palais Smale type. 

\begin{defin} 
A functional $J$ satisfies the PS condition at the level $m$ if for every sequence $u_n$ such that 
$$J(u_n) \to m, \quad dJ(u_n)\to 0$$
 the sequence $u_n$  has a converging subsequence.
\end{defin}

\begin{teo}[See {\cite{AM}}]\label{questo}
If $J$ is convex, it is of class $C^{1,1}$, satisfies (PS) condition and is bounded below, 
there exists $\lim_{t \to \infty} \tilde u = u$ and $u$ is a minimum of the functional $J$. 
In particular it satisfies $dJ(u)=0$. 
\end{teo}

These results can be directly applied to the functionals we have introduced, and  will be intensively used in the next sections.

\begin{ex}\label{ex1}
We consider the functional $J$ defined in \eqref{J} on the Sobolev space $W^{1, 2}_{0, \Lambda}$. 
Here we assume that $\theta$ is measurable, and $a_i$ are measurable and bounded. 
Then functional is of class $C^{1,1}$ and its differential~is
\begin{equation}\label{dj}
dJ(u)(h)=dF(u)(h)-  \int_Q  f h= \int_Q \sum_{i =1}^2  a^2_{i}  X_{\theta(x), i} u X_{\theta(x), i}h - \int f h.
\end{equation}
If both $a_{i}$ are different from $0$ the convergence of the steepest descent method is known. 
\end{ex}

If $a_1=1$,  $a_2 =0$ and 
$\theta(x,y)$ is a $C^{1}$ approximation of one of the functions 
defined in \eqref{pin} or \eqref{salt}, the operator is totally degenerate. 
 In this case we recover the operator 
 \eqref{sub} associated to V1 receptive profiles. 
Also, in this case it is possible to prove the (PS) condition. Indeed, if we consider a sequence such that 
$dJ(u_n)(h)\to 0$ and $J(u_n)\to \inf_{W^{1, 2}_{0, \Lambda}} J, $ then $ \norm{u_n}_{W^{1, 2}_{0, \Lambda}}^2$ is bounded, so that there is a function $u\in W^{1, 2}_{0, \Lambda}$ such that $  u_n$ weakly converge to  $  u$ in $W^{1, 2}_{0, \Lambda}$
and $u_n$ weakly converges to $u$ in $L^2$. From $dJ(u_n)\to 0$ we also deduce that 
$dJ(u_n)(u_n)\to 0$, so that 
 $$\lim_n \norm{u_n}_{W^{1, 2}_{0, \Lambda}}^2  = \lim_n  \int f u_n= \int f u,$$
 for the weak convergence. From the assumption of the differential we also know that $dJ(u_n)(u)\to 0$,  so that
  $$\norm{u}^2_{W^{1, 2}_{0, \Lambda}}  = \lim_n \int_Q X_{\theta(x), 1} u_n X_{\theta(x), 1}u  = \int f u.$$
It follows that $u_n \to u$ weakly in $W^{1, 2}_{0, \Lambda}$ and  $\norm{u_n}_{W^{1, 2}_{0, \Lambda}}^2 \to \norm{u}_{W^{1, 2}_{0, \Lambda}}^2$, which implies the convergence of $u_n$ to $u$ in 
$W^{1, 2}_{0, \Lambda}$.

As a consequence, by Theorem \ref{questo} the solutions of the stationary problems are the limit of the solution of the parabolic one. 

\begin{ex}
The operator in \eqref{general} is
the differential of the following functional 
$$ J(u)=\sum_{i=1}^2 \int _Q  a_i^2({X_{\theta(x), i}} h)^2(x)  dx +\int _Q  (\Delta u)^2(x)  dx - \int _Q f(x) u(x). $$
As before, this operator satisfies the (PS) condition. Indeed, a sequence such that $dJ(u_n)\to 0$ and $J(u_n)\to \inf J$ 
has the Laplacian bounded in $L^2$, so that, arguing as before, we obtain the convergence of $u_n$. The assumptions of Theorem \ref{questo} are satisfied, 
which ensures that the solution obtained through the parabolic approximation tends to the solution of the forth order operator. \newline
\end{ex}

We applied here the steepest descent method: the main limitation of the method is that it can be applied 
only to minima of functionals. However it has the big advantage that it can be applied to degenerate operators as the one of 
Example \ref{ex1}, since the functional, while approximating a minimum, bounds the norm of the minimizing sequence. 
Other powerful methods, as for example the semigroup technique, allow to handle the case where the equation does not come from a minimum of a functional. However the convergence result requires that the biggest eigenvalue of the operator is negative, which is not the case of the degenerate operators in Example \ref{ex1}.

\subsubsection{The Green function} 
\label{green}
The previous method allows to find the Green function of the operator $L_\Lambda$ an open, smooth sets $Q\subset \R^d$. The Green function is a smooth kernel 
$\Gamma\colon Q\times Q\backslash \{(x,x): x \in Q\}\to \R$ 
such that any solution of the problem 
$$L_{\Lambda} u =f, \text{  in  } Q, \quad  \partial_\nu u = 0 \text{  on  } \partial Q$$
can be represented in the form 
$$u(x) = - \int_{Q}\Gamma(x,y)f(y)dy.$$

For large class of 
operators with measurable coefficients analogous to $L_{\Lambda}$, the existence of a Green function, or fundamental solution (which is the analogous kernel on the whole space), is known. 
In \cite{GT} this is proved for second order uniformly elliptic operators, while the existence of a fundamental solution for uniformly sub-elliptic operators, is due to \cite{RS}.
 
If we fix the first entry $x\in Q$, the Green function as a function of its second entry is a smooth function  $\Gamma(x, -)\colon Q\backslash \{x\}\to R$ solution of the problem
$$L^*_\Lambda \Gamma(x, -)\ = \delta_x \text{ in } Q, \quad 
\partial_\nu \Gamma=0\text{ on }\partial Q,$$
where $L^*_\Lambda$ is the formal adjoint of $L_\Lambda.$

In order to solve this last problem and find the Green function, we can apply the steepest descent method we have introduced (see also results in Section \ref{NR}) .

\section{Discretization and homogenization results}\label{DC}

Scope of this section is to provide a homogenization results for operators of the same type as the one defined in \eqref{operats} and \eqref{op}. 
We can express the theory in an open set $Q \subset \R^d$. To this end we will consider a finite set $\Lambda\subset \Z^d$ symmetric with
respect to~$0$, and we will call $m$ the number of its elements. Any element $z_i$ with $i=1,\cdots,  m$ is a vector,   
the associated directional derivative will be denoted as $X^\varepsilon_{z_i}$. 
We intend to study the discretization of the following Dirichlet problem
which generalizes the operator \eqref{op} defined above:
\begin{equation}\label{contD}
-\sum_{i=1}^m(X^\varepsilon_{z_i})^*\left(a_{ij}(x) X^\varepsilon_{z_j}  u(x)\right)=f(x).
\end{equation}

We will consider both  deterministic and stochastic coefficients $a_{ij}$.  
For deterministic coefficients we study the convergence of the discretized problem to the continous one. 
In the stochastic case, we show that, even if we start with a random orientation 
$\theta$ in the discrete problem, a homogenization procedure will ensure the convergence to a 
continuous problem with constant coefficients.

\subsection{Difference operators}

We recall here how to discretize this second order differential operator 
referring to \cite{koz87} and \cite{piatremy01} for a wider introduction.

We start with 
discretizing a smooth bounded domain $Q\subset \R^d$. 
Let us call $Q_\varepsilon=Q\cap \varepsilon\Z^d$, where $\varepsilon>0$.
The boundary $\partial Q_\varepsilon^\Lambda$ of $Q_\varepsilon$ is defined by
$$\partial Q_\varepsilon^\Lambda:=\{x+\varepsilon z_i|\ x\in Q_\varepsilon,\ i=1, \cdots m\}\backslash Q_\varepsilon.$$
We also introduce $\overline Q_\varepsilon:=Q_\varepsilon \cup \partial Q_\varepsilon^\Lambda$.

Consider a function $u^\varepsilon\colon Q_\varepsilon\to \R$, then the standard
difference operator associated $X_{z_i}^\varepsilon$ is defined by
$$X_{z_i}^\varepsilon u^\varepsilon(x):=\frac{u^\varepsilon(x+\varepsilon z_i)-u^\varepsilon(x)}{\varepsilon}\ \ \forall  i=1, \cdots m.$$
We explicitly note that $z_i$ is a vector, not necessary an element of the canonical basis, 
of the space, so that in general $X_{z_i}^\varepsilon $ represent the discretization of the directional derivative $X_i$, not a partial derivative.

Here we can consider a matrix $A^\varepsilon(x)=(a^\varepsilon_{ij}(x))_{ij=1, \cdots m}$ 
when $x\in Q_\varepsilon$. With these notations, the discrete Laplace operator can be written as follows
\begin{equation}\label{dir}
L^\varepsilon_\Lambda u^\varepsilon(x):=\sum_{i=1}^mX^\varepsilon_{-z_i}(a^\varepsilon_{i j}(x)X^\varepsilon_{z_i}u^\varepsilon(x)).
\end{equation}
Consequently, the 
Dirichlet problem 
with second member $f^\varepsilon\colon Q_\varepsilon \to \R$ becomes 
\begin{equation}\label{dir}
L^\varepsilon_\Lambda u^\varepsilon(x)=f^\varepsilon(x)\ \ \forall x\in Q_\varepsilon,
\end{equation}
with $u^\varepsilon(x)=0\ \m{if}\ x\in \partial Q_\varepsilon^\Lambda.$

\begin{defin}
We say that the problem (\ref{dir}) is uniformly elliptic if there exist constant $c_1,c_2,\varepsilon_0>0$ such
that for any $\eta \in \R^d$ represented as $\eta = \sum_{i=1}^m \eta_i z_i$ one has 
\begin{align}
|a^\varepsilon_{i j}(x)|&\leq c_1\ \ \forall x\in \overline Q_{\varepsilon}\ \forall i,j =1, \cdots m\\
c_2\cdot \norm {\eta}^2_E &\leq \sum_{i,j=1}^m
a^\varepsilon_{i,j}(x)\eta_i \eta_j \; \;  \forall x\in \overline Q_{\varepsilon}
\end{align}
where $\norm {\eta}_E$ is the Euclidean norm. 
\end{defin}

Let us now define suitable functional spaces where we will look for the solution. 
We consider $v^\varepsilon\colon \varepsilon \Z^d\to \R$, 
and we define the analog of the $L^2$-norm
$$\norm{u^\varepsilon}^2_{L^2(Q_\varepsilon)}:=\varepsilon^d \ \sum_{x\in Q_\varepsilon}|v^\varepsilon(x)|^2.$$

If we consider $v^\varepsilon\colon \varepsilon \Z^d\to \R$, we say that $v^\varepsilon$ is in $W^{1,2}_0(Q_\varepsilon)$
if $v^\varepsilon(x)=0$ for any $x\notin Q_\varepsilon$, this is in fact a translation of the
notion of Sobolev space to the discrete setting. 
The notion of discrete gradient descends from this
setting straightforwardly.

\begin{defin}
Consider $v^\varepsilon \in W^{1,2}_0(Q_\varepsilon)$,
$$\nabla^\varepsilon_\Lambda v^\varepsilon(x):= \sum_{j=1}^m X^\varepsilon_{z_j} v^\varepsilon(x)X^\varepsilon_{z_j}\ \ \forall x\in\overline Q_\varepsilon.$$
\end{defin}

Due to the uniformly ellipticity condition, we can choose the following norm on $W^{1,2}_0$,
$$\norm{v^\varepsilon}^2_{W^{1,2}_0(Q_\varepsilon)}:=\varepsilon^d\cdot \sum_{x\in \overline Q_\varepsilon}\sum_{i=1}^m,
\left|X^\varepsilon_{z_j} v^\varepsilon(x)\right|^2.$$
and denote by $W^{-1,2}(Q_\varepsilon)$ the dual space to $W^{1,2}_0(Q_\varepsilon)$.\newline


We will provide a definition of solutions of the equation of the operator $L^\varepsilon_\Lambda$ generalizing to the discrete setting a mean value formula. 
Indeed it is well known that if $\Delta$ is the standard 
Laplacian, $B(x,\varepsilon)$ is the ball of radious $\varepsilon$ and 
$|B(x,\varepsilon)|$ is its Lebesgue measure, 
then 
$$u(x) - \frac{1}{|B(x,\varepsilon)|}\int_{B(x,\varepsilon)} u(y) dy = 
\frac{\varepsilon^2}{4 + 2d} \Delta u + o(\varepsilon^2)
\text{ as } \varepsilon \to 0.$$
The same formula still holds for large classes of linear operators.
%
%
%

In the discrete setting we will consider the following formula introduced in \cite[Proposition 1.3]{piatremy01}.

\begin{prop}\label{prop1}
Consider a function $p_{z_i}^\varepsilon\colon Q_\varepsilon\to \R$ defined for
any $z_i\in \Lambda$, if it satisfies the following three properties
\begin{enumerate}
\item $p_{z_i}^\varepsilon(x)\geq 0$ and $\sum_{i=1}^mp_{z_i}^\varepsilon(x)=1\ \forall x\in Q_\varepsilon$;
\item $\exists \delta>0$ such that $p_{\pm e_i}^\varepsilon\geq \delta$ for $i=1,\dots,d$;
\item $p_{z_i}^\varepsilon(x)=p_{-z_i}^\varepsilon(x+\varepsilon z_i)$;
\end{enumerate}
then the following problem is uniformly elliptic
\begin{equation}\label{eq_pe}
u^\varepsilon(x)-\sum_{i=1}^mp_{z_i}^\varepsilon(x)u^\varepsilon(x+\varepsilon z_i)=\varepsilon^2\cdot f^\varepsilon(x)
\ \ \m{in}\ Q_\varepsilon,\ \ \ u^\varepsilon(x)=0\ \forall x\in \partial Q_\varepsilon^H,
\end{equation}
for any $f^\varepsilon\colon Q_\varepsilon\to \R$.
\end{prop}
\begin{rmk}
The operator associated to the previous mean value formula can be rewritten in the form (\ref{dir})
with $a^\varepsilon_{ij}(x)=p_{z_i}^\varepsilon(x)$ if $z_i=z_j\neq0$ and $0$ otherwise.\newline
\end{rmk}

In order to compare the functions defined over $Q_\varepsilon$ with those
having a continuous argument, we follow Kozlov \cite{koz87} in defining a mesh completion.
In particular if $f^\varepsilon \colon Q_\varepsilon \to \R$ then we denote by $\tilde f^\varepsilon\colon Q\to \R$
the function such that
$$\tilde f^\varepsilon(x)=f(y_x)\ \ \forall x\in Q,$$
where $y_x$ is the  point in $Q_\varepsilon$ of components $y_x =(y_1, \cdots, y_d)$ such that
$$y_i-\frac{\varepsilon}{2}\leq x_i< y_i+\frac{\varepsilon}{2}\ \ \forall i=1,\dots, d.$$
In what follows, if it is clear from the context we will use 
a little abuse of notation by denoting $\tilde f^\varepsilon\in W^{1,2}_0(Q)$ also by
$f^\varepsilon$. Furthermore, we will say that $f^\varepsilon$ converges (strongly or weakly)
to $f$ in $L^2(Q_\varepsilon)$, $W^{1,2}(Q_\varepsilon)$ or $W^{-1,2}(Q_\varepsilon)$
when the mesh completion $\tilde f^\varepsilon$ converges to $f$ in $L^2(Q)$,
$W^{-1,2}(Q)$.
 If needed we will replace this mesh completion with a piecewise linear one $\tilde f^\varepsilon$ in $W^{1,2}(Q)$. In this case we will say that  $f^\varepsilon$ converges (strongly or weakly)
to $f$ in $W^{1,2}(Q_\varepsilon)$ if the piecewise linear mesh completion $\tilde f^\varepsilon$ converges to $f$ in 
$W^{1,2}(Q)$.

%


\subsection{Convergence results - deterministic model}
We are interested in treating the convergence of difference operators
to  usual differential operators
when imposing $\varepsilon\to 0$. 
We start
by considering an existence result and its uniform estimates.
\begin{prop}
If the problem (\ref{dir}) is uniformly elliptic and $f^\varepsilon\in L^2(Q_\varepsilon)$,
then there exists unique a solution of the problem $u^\varepsilon\in W^{1,2}_0(Q_\varepsilon)$
and there exists $c>0$ such that
$$\norm{u^\varepsilon}_{W^{1,2}_0}\leq c\cdot \norm{f^\varepsilon}_{L^2},$$
uniformly in $\varepsilon$.
\end{prop}

In what follows we denote as usual by $\rightharpoonup$ the weak convergence.
The weak convergence of difference operators preserves in fact the 
notion of partial derivative.

\begin{prop}[{\cite[p.355]{koz87}}]
If $u^0\in W^{1,2}_0(Q)$ (or $L^2(Q)$) and $u^\varepsilon\rightharpoonup u^0$ in $W^{1,2}_0(Q_\varepsilon)$
(or $L^2(Q_\varepsilon)$), then
$$X^\varepsilon_{e_i}u^\varepsilon\rightharpoonup \frac{\partial u^0}{\partial e_i}\ \ \m{in}\ L^2(Q_\varepsilon)\ \m{or}\ W^{-1,2}(Q_\varepsilon).$$\newline
\end{prop}

We will denote by $A$ an uniformly elliptic matrix and by 
 $A^\varepsilon$ its discretization on the $\varepsilon$ 
grid.
Consider a family of uniformly elliptic discrete Dirichlet problems
(\ref{dir}). 
As before we denote by $A^\varepsilon(x)$ 
the coefficient matrix. We will also denote by $A(x)$ a $m \times  m$ real matrix 
defined for $x\in Q$, whose coefficients
will be denoted $a_{ij}$.\newline

Denote by $u^\varepsilon$ the solution of the Dirichlet problem
associated to $u^\varepsilon$ and by $u^0$
the solution of the continuous Dirichlet problem
\eqref{contD}.
\begin{defin}\label{410}
We say that the matrix $A^\varepsilon$ $H$-converges to $A$,
$A^\varepsilon\xrightarrow[\varepsilon\to 0]{H}A$,
if for any sequence $f^\varepsilon \in W^{-1,2}(Q_\varepsilon)$
such that $f^\varepsilon\to f\in W^{-1,2}(Q)$, we have 
\begin{align*}
u^\varepsilon &\rightharpoonup u^0 {\text {  in }} W^{1,2}_0(Q_\varepsilon) \\
D_{a^\varepsilon} u^\varepsilon =\sum_{j=1, \cdots m}a^\varepsilon_{i j }X^\varepsilon_{j}u^\varepsilon&\rightharpoonup
D_{a} u^0 =\sum_{j=1}^ma_{ij }X_j u^0  {\text {  in }} L^2(Q_\varepsilon)
\end{align*}
\end{defin}

\begin{defin}
Given a matrix valued function $A^1(x)=(a^1_{ij}(x))$ with $i,j =1, \cdots m$ and $x\in \Z^d$,
and a sequence of induced matrices
 $A^\varepsilon(x):=A^1(x\slash \varepsilon)$ for any $x\in Q_\varepsilon$,
if the associated problem (\ref{dir}) is uniformly elliptic and 
there is a constant matrix $A^0$ such that 
$A^\varepsilon\xrightarrow{H} A^0$,
then we call the matrix $A^0$ the homogenized matrix for $A^\varepsilon$.
\end{defin}

\subsection{Convergence results - stochastic model}

We will now consider an operator with random coefficients, similar to the one introduced in 
\eqref{salt} and in Section \ref{secordop} to describe the action of the receptive profiles.

In order to treat the case of random coefficients, consider
a probability space $(\Omega, \mc F, \mu)$ and a group
 $\{T_x\colon \Omega\to\Omega|\ x\in \Z^d\}$
of $\mc F$-measurable transformations
respecting the following properties:
\begin{enumerate}
\item $T_x\colon \Omega\to \Omega$ is $\mc F$-measurable $\forall x\in \Z^d$;
\item $\mu(T_x\mc B)=\mu(\mc B)$, $\forall \mc B\in \mc F$ and $\forall x\in \Z^d$;
\item $T_0=\id$ and $T_x\circ T_y=T_{x+y}\ \forall x,y\in \Z^d$.
\end{enumerate}
\begin{defin}
The group $T_x$ is ergodic if for any $f\in L^1(\Omega)$
such that $f(T_x\omega)=f(\omega)$ for $\mu$-a.e.~$\omega\in \Omega$ and for any $x\in \Z^d$, there exists
a constant $K$ such that $\mu$-a.s.~$f=K$.
\end{defin}

In what follows we suppose the group $T_x$ to be ergodic
and we build a family of random operators. Consider a
 $\mc F$-measurable function with matrix values
$\mc A(\omega)=(a_{ij}(\omega))$ for $\omega\in \Omega$ such that $\mc A$ is a $m\times m$ symmetric matrix.
We define the family of operators
$$A^\varepsilon(x)(\omega):=\mc A(T_{x\slash \varepsilon}\omega),\ \ \forall  \omega\in \Omega,  x\in \varepsilon\Z^d.$$

\begin{rmk}\label{rmk_A}
We state a condition on $\mc A$ that implies the uniform ellipticity of the family~$A^\varepsilon$.
Consider $\eta \in \R^d$ represented as $\eta= \sum_{i=1}^ m \eta_i z_i$.
If there exist $c_1,c_2>0$ such that the following inequalities
are true a.s.
\begin{align*}
 |a_{ij}(\omega)|&\leq c_1, \\
c_2 \cdot \norm{\eta}_{E}&\leq 
\sum_{z_i,z_j\in \Lambda\backslash \{0\}}
a^\varepsilon_{ij}(\omega)\eta_i\eta_j,\ \forall x\in Q_\varepsilon
\end{align*}
for all $\eta \in \R^d$ represented as $\eta= \sum_{i=1}^ m \eta_i z_i$
then the $A^\varepsilon$ are uniformly elliptic in any regular domain $\Omega$.
\end{rmk}




\begin{teo}[{\cite[Theorem 2.17]{piatremy01}}]
If the operators $A^\varepsilon$ are built as above, with $T_x$ ergodic group
and the matrix $\mc A$ respecting the conditions of Remark \ref{rmk_A},
then a.s. the family $A^\varepsilon$ admits a homogenization and the homogenized
matrix $A^0$ does not depend on $\omega$.\newline
\end{teo}

 Let us now apply our theoretical results to equation \eqref{salt}, and to the choice of
  orientation \eqref{pin} under the simplified assumption that the only possible direction are vertical or horizontal.

We consider a particular family of operators $A^\varepsilon$ and analyze 
its asymptotic behavior. In particular we work in dimension $d=2$
and use some results of percolation theory
in order to show that for some well defined families,
asymptotically the operator converge to a Laplacian.\newline

{\bf Randomly chosen horizontal or vertical orientation}

We split $\R^2$ into squares $\left\{\left[-\frac{1}{2},\frac{1}{2}\right]+j|\ j\in \Z^2\right\}$,
and consider a random variable $\kappa_r$ defined over $\R^2$,
such that if $\delta>0$, then
$$\kappa_r=\left\{\begin{array}{l}
\delta\ \ \text{with probability }r\\
1\ \ \text{with probability }1-r
\end{array}\right.$$
where
$0<r<1$. Given a finite subset $\Lambda\subset \Z^2$ symmetric with respect
to $0$, we consider the transition functions $p_{z_i}(x)$ for any $x\in \R^2$ and $z_i\in \Lambda$,
and suppose the $p_{z_i}(x)$ are functions defined on the set $\{\kappa_r(x+z_i)|\ z_i\in\Lambda\}$.

The independece of $\kappa_r(j)$ for different $j\in \Z^2$ implies that 
the transformation group allowing the construction of the families $\{p_{z_i}(x)|\ x\in \R^2\}$
is ergodic, or equivalently that if $\forall z_i\in \Lambda$ if $p_{z_i}(x)=p_{z_i}(0)$ a.s., then a.s. $p_{z_i}$ equals a constant.

From the $p_{z_i}(x)$ we define the transition functions $p^\varepsilon_{z_i}(x)$
in order to describe a uniformly elliptic problem as in the equation (\ref{eq_pe}),
$$p^\varepsilon_{z_i}(x):=p_{z_i}(x\slash\varepsilon)\ \ \forall x\in Q_\varepsilon\subset \varepsilon\Z^2.$$
The operators $A^\varepsilon=(a^\varepsilon_{ij})$ are defined by 
$a^\varepsilon_{ii}(x)=p^\varepsilon_{z_i}(x)$ if $z_i\neq 0$ and $a^\varepsilon_{ij}(x)=0$ otherwise.

\begin{prop}[{\cite[\S2]{koz87}}]
If the functions $p_{z_i}^\varepsilon$ built as above satisfy a.s.~the three
hypothesis of Proposition \ref{prop1}, then
the $A^\varepsilon$ $H$-converge a.s. to an elliptic
operator $A^0$ with constant coefficients and moreover
$A^0$ is isotropic,
$$A^0=(a^0_{ij})=a^\delta(r)\cdot \id,$$
where $a^\delta(r)\in \R$.\newline
\end{prop}

The work \cite{piatremy01} focus on the case 
where 
 the functions
$p_{z_i}$ are defined by

$$p_{z_i}(x):=\left\{\begin{array}{cl}
\frac{2\kappa_r(x)\cdot\kappa_r(x+z_i)}{4(\kappa_r(x)+\kappa_r(x+z_i))} & \m{ if } z_i\in \Lambda\backslash\{0\}\\
1-\sum_{z_i\neq 0}p_{z_i}(x) & \m{ if }z_i=0.
\end{array}\right.$$
%

\begin{figure}[htbp]
\label{f4x}
\centering
\includegraphics[height = 1.7in ]{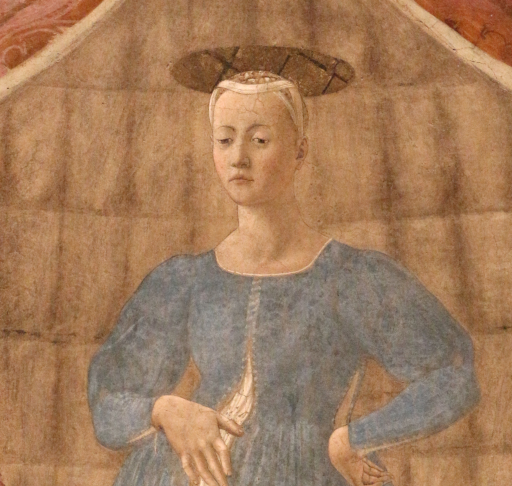}  
\includegraphics[height= 1.7in ]{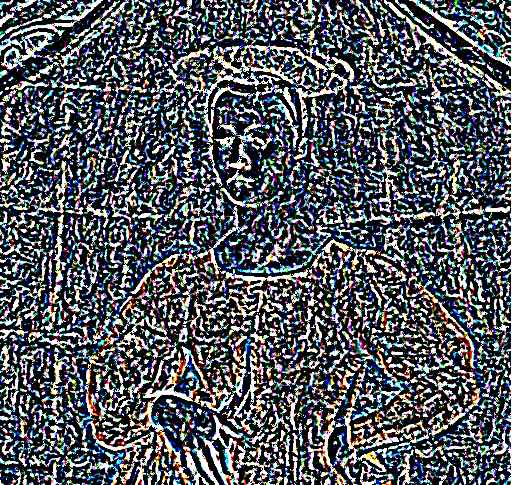}
\caption{Original painting of Piero della Francesca (left) and its second order differentiaton with randomly chosen horizontal or vertical orientation (right).\newline} 
\label{f4}
\end{figure}

{\bf Second and higher order operators}

Let us explicitly observe that the convergence theory developed up to here can be applied 
only to second order operators. It is not clear if similar results could be 
obtained for operators of order different from a point to the other.

\section{Results} \label{NR}

The reconstruction model  of Equation (\ref{maineq}) is implemented by  numerically solving the parabolic Equation (\ref{claprob1})

\begin{equation}\label{resultheat}
u_t= L_\Lambda u -  L_\Lambda I
\end{equation}

with the heterogeneous operator $L_\Lambda$ given by Equation  \eqref{operats} and  $\beta = 2$, so that to take into account three main kind of receptive profiles:

\begin{equation}\label{resulteq}
L_\Lambda = a_1\Delta+ a_2 X^2_{\theta,1} - a_3 X^4_{\theta,1}.
\end{equation}

The spatial domain $Q$ is discretized with uniform spacing $\varepsilon$ to obtain the matrix $Q_\varepsilon$ on which the solution $u^\varepsilon$ is defined. 
 Centered finite differences are used to approximate spatial derivatives of the Laplacian 
 $$\Delta=\partial^2_{x_1x_1} + \partial^2_{x_2x_2} $$
as well as the second orderd sub-Riemannian term  
$$X^2_{\theta,1} = \cos^2(\theta)\partial^2_{x_1x_1} + 2 \cos(\theta)\sin(\theta)  \partial^2_{x_1 x_2} + \sin^2(\theta)\partial^2_{x_2x_2},$$
where $(x_1, x_2)$ denotes a point in $Q_\varepsilon$. 
The fourth order sub-Riemannian term is computed by applying twice the second order term,
$$X^4_{\theta,1}=X^2_{\theta,1}(X^2_{\theta,1}).$$

Time evolution is discretized with forward finite difference and Neumann boundary conditions are applied.
Typical time step of the evolution is $dt=0.1$ if just second order operators are involved, while $dt=0.001$ is adopted if fourth order operators are present,
in order to garantee the stability of finite difference method \cite{strik89}. The stopping criterion is directly related to the error of convergence of the solution so that the evolution stops when $\sum_{x\in Q_\varepsilon}|u^\varepsilon_{k+1}(x)-u^\varepsilon_k(x)|<\varepsilon_c$, where $u^\varepsilon_{k}$ and $u^\varepsilon_{k+1}$ are solutions corresponding to two subsequent time steps and $\varepsilon_c=10^{-4}$.

The three bands of the RGB image stimulus shown in Figure \ref{f4} left have been processed separately. 
Different mixtures of second and fourth order operators in Euclidean and sub-Riemannian metrics have been considered by changing the parameters $a_1, a_2, a_3$.\newline

\begin{figure}[htbp]
\label{f5x}
\centering
\includegraphics[height = 1.5in ]{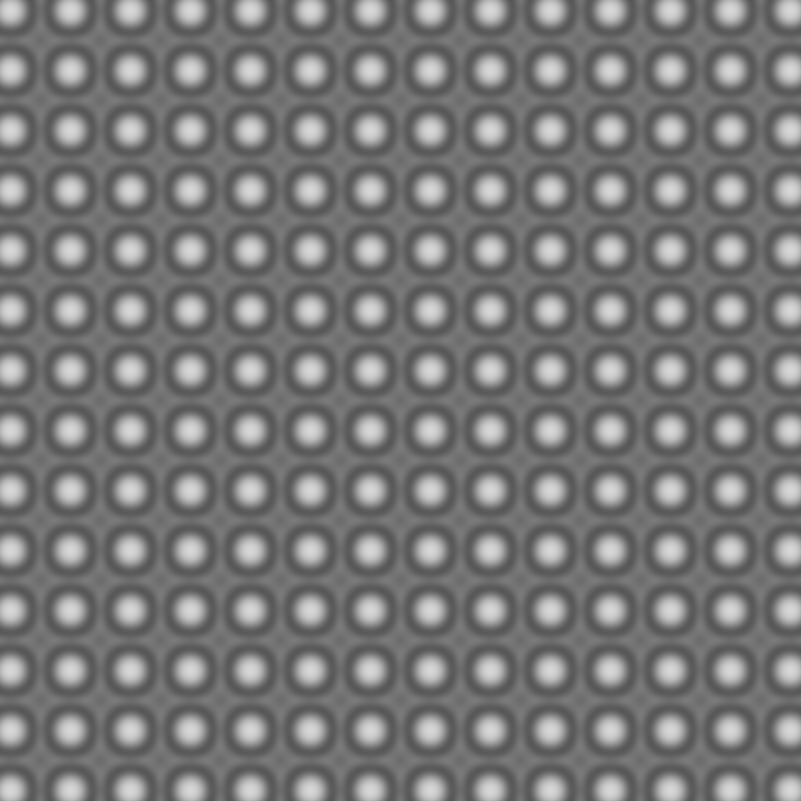}
\includegraphics[height = 1.5in ]{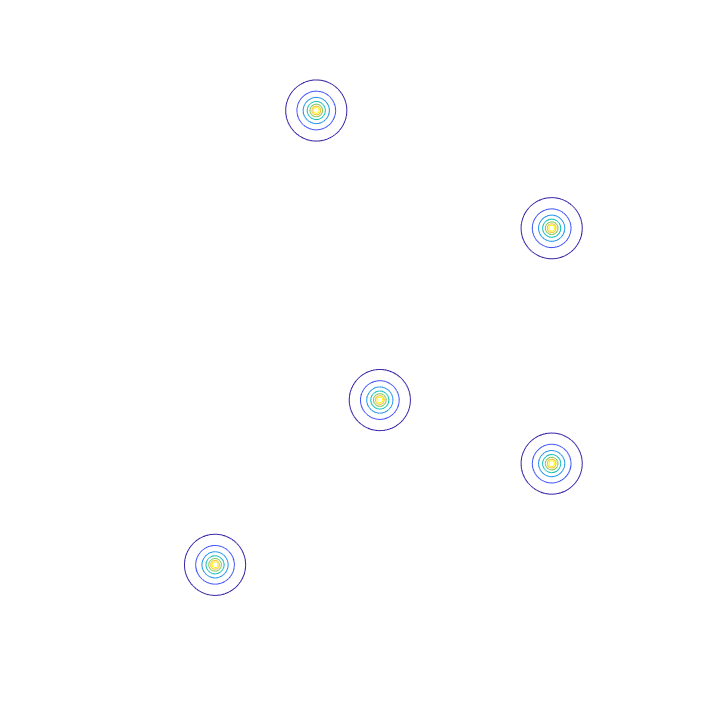} \\
\includegraphics[height = 1.4in ]{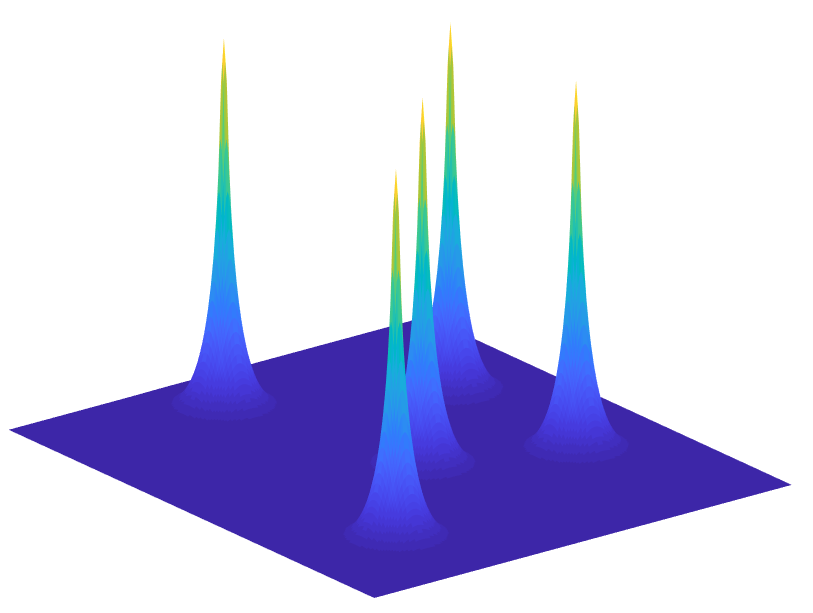}
\includegraphics[height = 1.5in ]{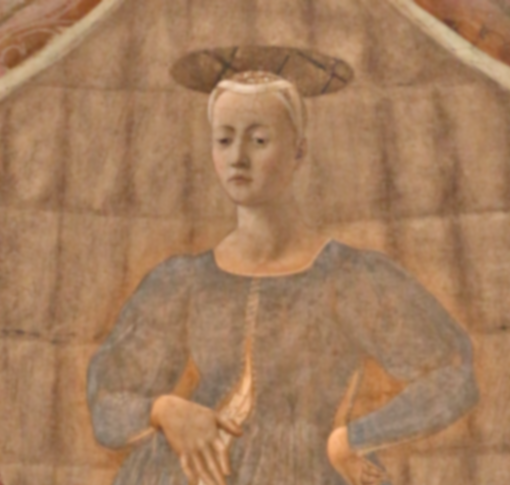}
\caption{Classical Laplacians. Upper left:  Ps are visualized in a subsampled set of points. Upper right: Level lines of the fundamental solution sampled in some points. Bottom left: Surface of the fundamental solution. Bottom right: Reconstructed image.}
\label{f5} 
\end{figure}

In our first numerical test, just isotropic Laplacians have been considered (Figure \ref{f5}) setting $a_1=1, a_2=0, a_3=0$). 
The set of  Ps are visualized in a sub-sampled set of points (upper left) and they are all equal. Equation  \ref{maineq}  corresponds in this case to the classical Poisson problem introduced by J.~M.~Morel to solve the Retinex algorithm with differential instruments \cite{Mo11}. The solution is visualised (bottom right) and corresponds to the image stimulus up to a harmonic function. The Green function $\Gamma(x_1,x_2)$ (see Section \ref{green})  is computed for some points $\gamma$. Its level lines are visualized (upper right) as well as its surface (bottom left). 
In case of classic Laplacians the  Green function in 2D is analytically  known $\Gamma(r)= \log(r)$, where $r=\sqrt{x_1^2+x_2^2}.$

\begin{figure}[htbp]

\label{f7x}
\centering
\includegraphics[height = 1.5in ]{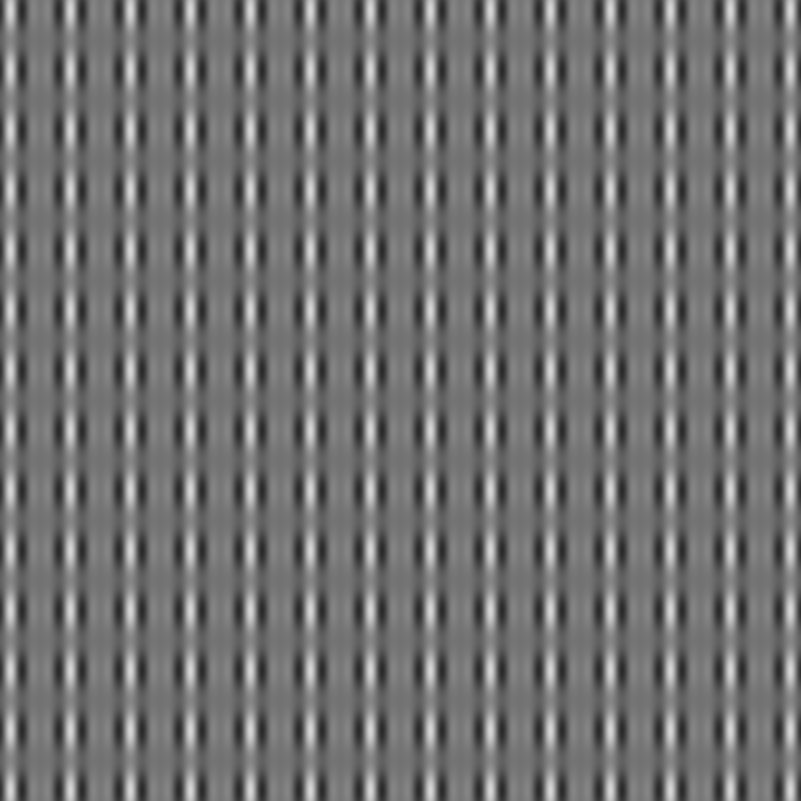} 
\includegraphics[height = 1.5in ]{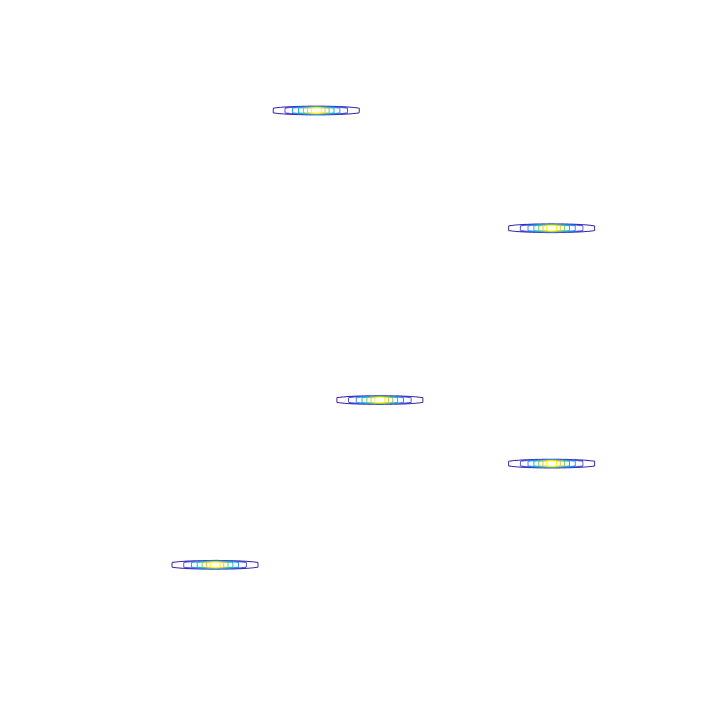}  \\
\includegraphics[height = 1.4in ]{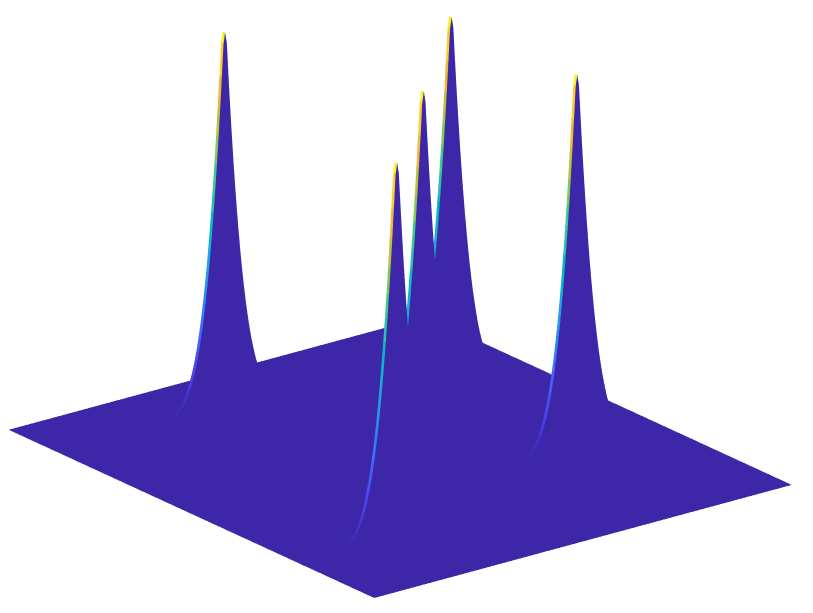}
\includegraphics[height = 1.5in ]{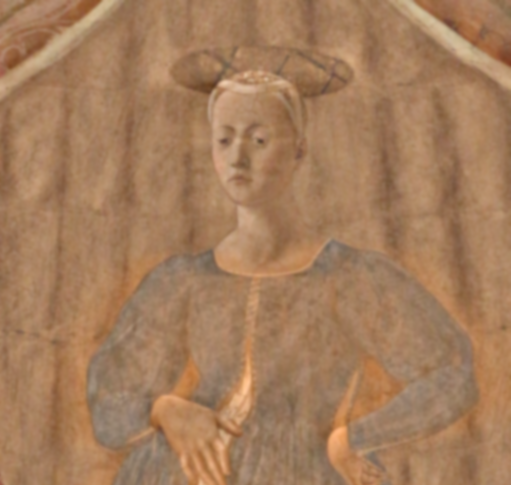}
\caption{Homogeneous sub-Riemannian second order operators $L_{\Lambda}=\partial^2_{x_1x_1}$. 
Upper left:  Ps are visualized in a subsampled set of points. Upper right: Level lines of the fundamental solution. Bottom left: Surface of the fundamental solution. Notice that fundamental solutions are strongly asymmetric. Bottom right: Reconstructed image:  Artifacts appear as uncorrelated horizontal lines.}
\label{f7}
\end{figure}

 If $a_1=0, a_2=1,a_3=0$ and $\theta=0$ in every point,  then
 the set of operators is constituted by directional derivatives with all the same orientation (Figure \ref{f7}) and  the corresponding parabolic equation is
 $$u_t= \partial^2_{x_1x_1} u -  \partial^2_{x_1x_1} I.$$  The Green function is very anisotropic and the reconstruction is corrupted by artifacts. This is due to the fact that many one dimensional Poisson problems are solved independently in the direction of the derivations and solutions are uncorrelated. 

\begin{figure}[htbp]
\label{f8x}
\centering
\includegraphics[height = 1.5in ]{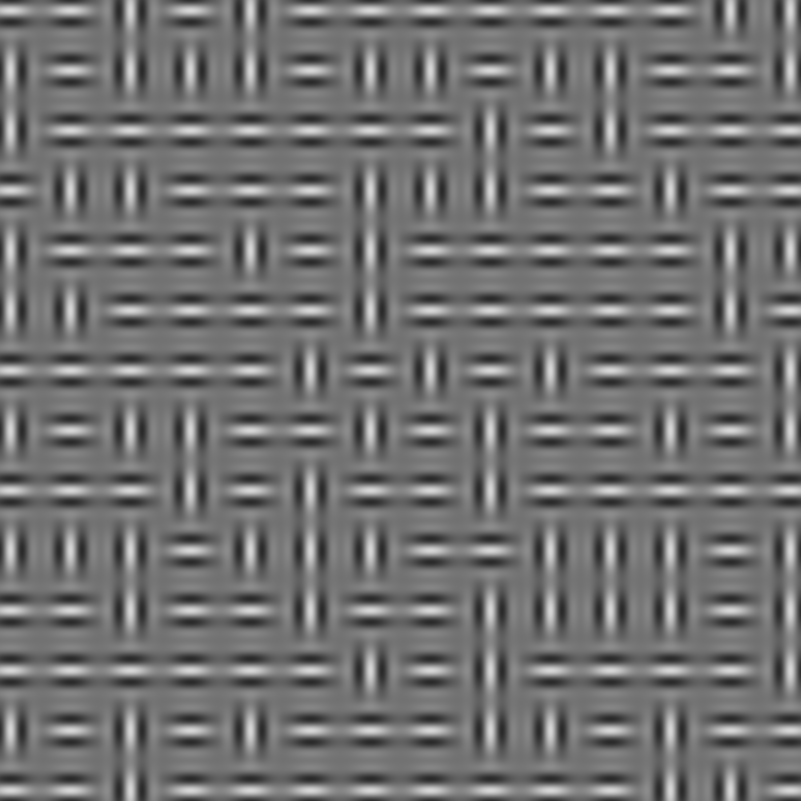}
\includegraphics[height = 1.5in ]{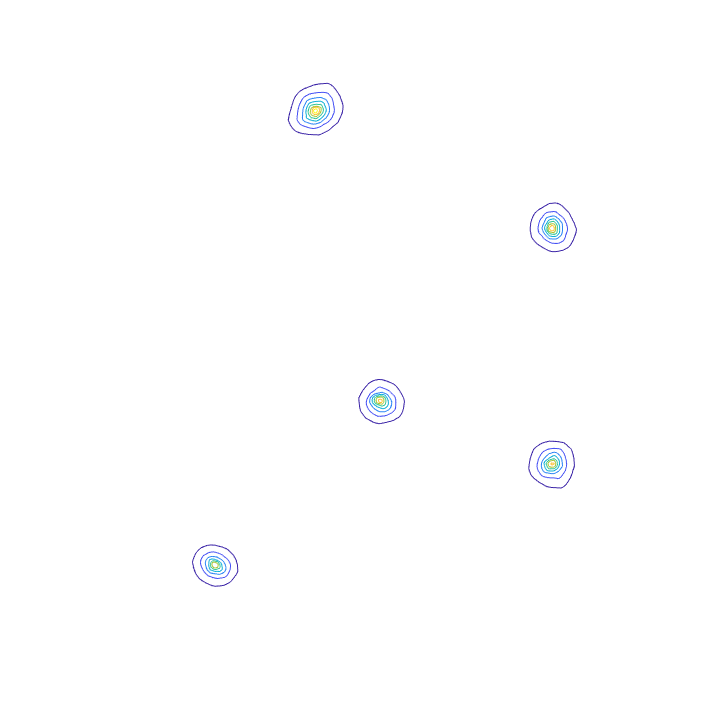}  \\
\includegraphics[height = 1.4in ]{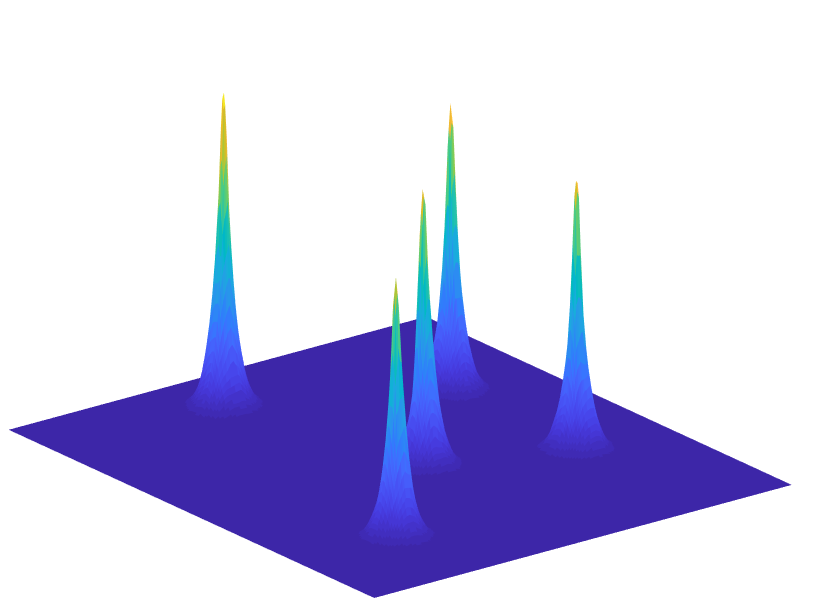}
\includegraphics[height = 1.5in ]{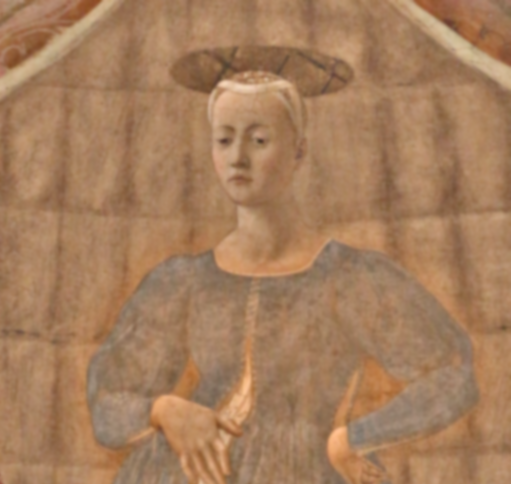} 
\caption{Sub-Riemannian second order operators in the vertical direction and horizontal directions. 
Upper left:  Ps are visualized in a subsampled set of points. Upper right: Level lines of the fundamental solution. Bottom left: Surface of the fundamental solution. Notice that fundamental solutions are different point by point and they are quite round  although they are composed by anisotropic operators. Bottom right: The reconstructed image appears free from artifacts.\newline} 
\label{f8}
\end{figure}

The case with a mixture of directional derivatives changing from point to point in the horizontal and vertical directions, considered in
 Figure \ref{f8}, is very different. 
As in the previous case,  $a_1=0, a_2=1,a_3=0$,  but now direction $\theta$ is randomly chosen 
depending on the position $x= (x_1,x_2)$ with values $\theta=0$ or $\theta=\pi\slash2$. 
The result of the differentiation of the stimulus image is shown in Figure \ref{f4} right. Notice that the Green functions of the mixture of operators change from
 point to point. Although the direction is randomly chosen, the Green function maintains a certain regularity and symmetry. 
 This allows to obtain a reconstruction of the perceived image without visible artifacts. 
 The formal proof of convergence of this distribution of operators towards classical Laplacians is given in Section~\ref{DC}.
 
\begin{figure}[htbp]
\label{f9x}
\centering

\includegraphics[height = 1.4in ]{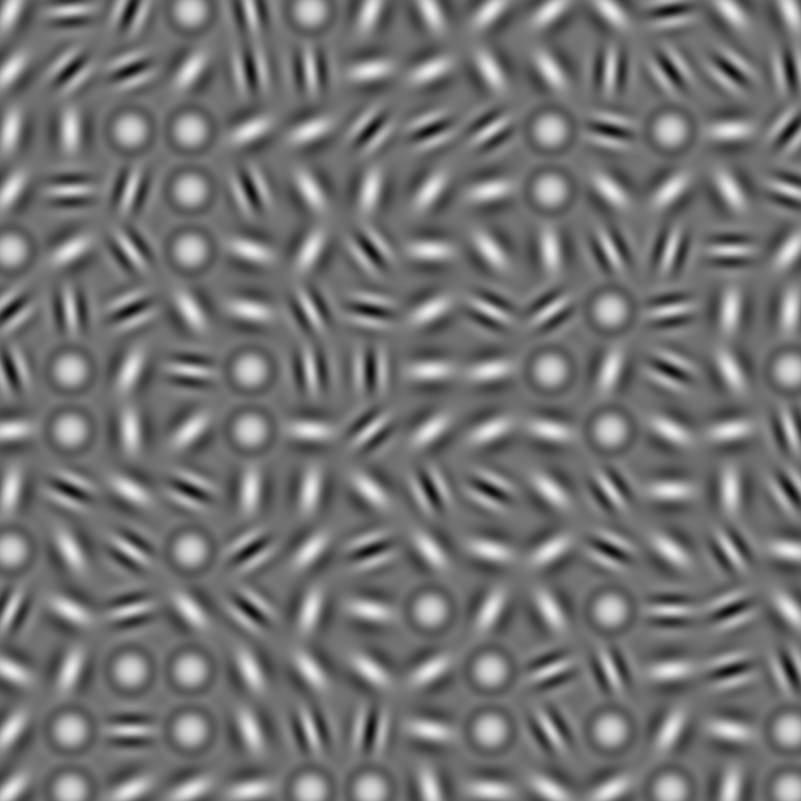}
\includegraphics[height = 1.4in ]{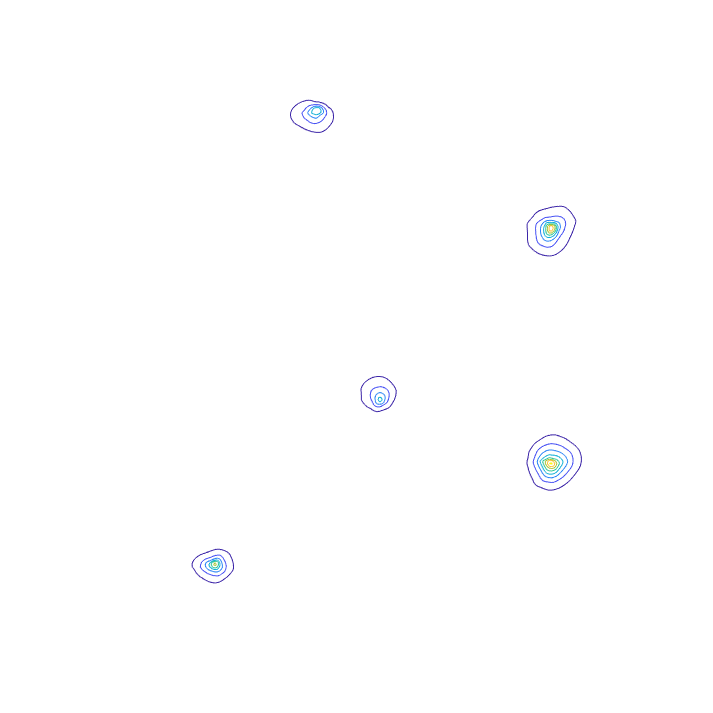}  \\
\includegraphics[height = 1.3in ]{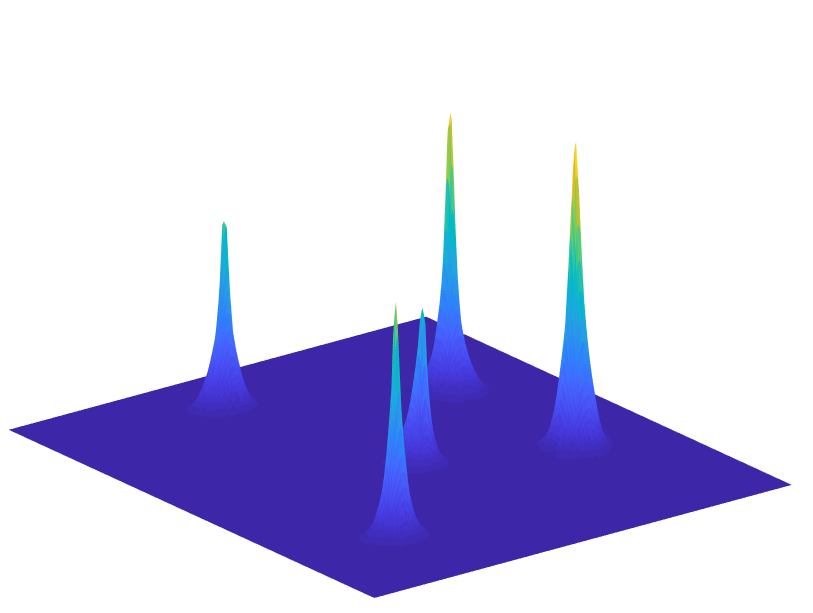}
\includegraphics[height = 1.4in ]{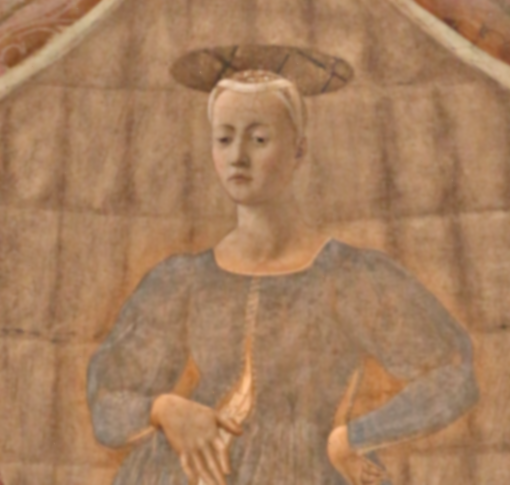}
\caption{Heterogeneous operators in the cortex of rodents: Either second or fourth order, Euclidean or sub-Riemannian operators are randomly distributed with random n-orientations. 
Upper left:  Ps are visualized. Upper right: Level lines of the fundamental solution. Bottom left: Surface of the fundamental solution. Notice that fundamental solutions are different point by point and they are quite symmetric although they are composed by very anisotropic operators of different degree. Bottom right: The reconstructed image appears free from artifacts.\newline} 
\label{f9}
\end{figure}

A typical distribution of Ps in rodents is a random salt and pepper mixture 
of Mexican hats as well as simple cells with different number of cycles (Figure \ref{f8}). 
These correspond to the case where the coefficients $a_i$ are chosen as a partition of the unit
such that
 in any point $(x_1,x_2)$ of the domain $Q_\epsilon$ just one parameter has value $1$ while the others are null. 
 The choice of the non-null parameter in every point is randomic, The distribution of the $\theta$ angles is also randomic over the spatial domain. 

Green functions of the mixture are different from point by point and they are quite symmetric although they are composed by very anisotropic operators of different degree. In this case the reconstructed image is the original stimulus $I(x_1,x_2)$ up to a harmonic function, where harmonic means the annulation of the heterogeneous operator $L_{\Lambda}=0$. 
No artifact is visible in the reconstructed-perceived image.

\begin{figure}[htbp]
\label{f10x}
\centering
\includegraphics[height = 1.4in ]{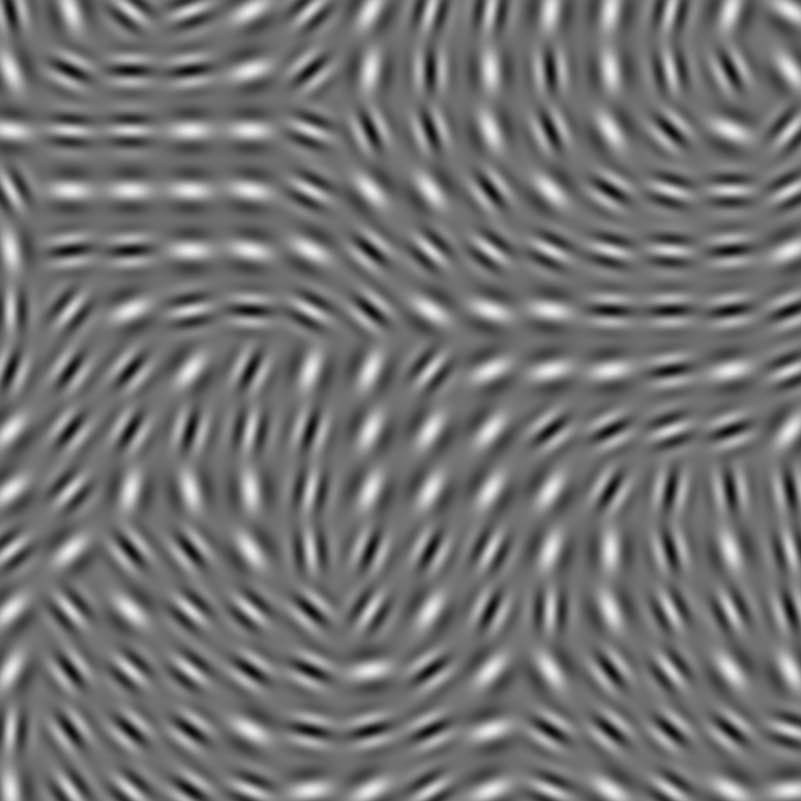}
\includegraphics[height = 1.4in ]{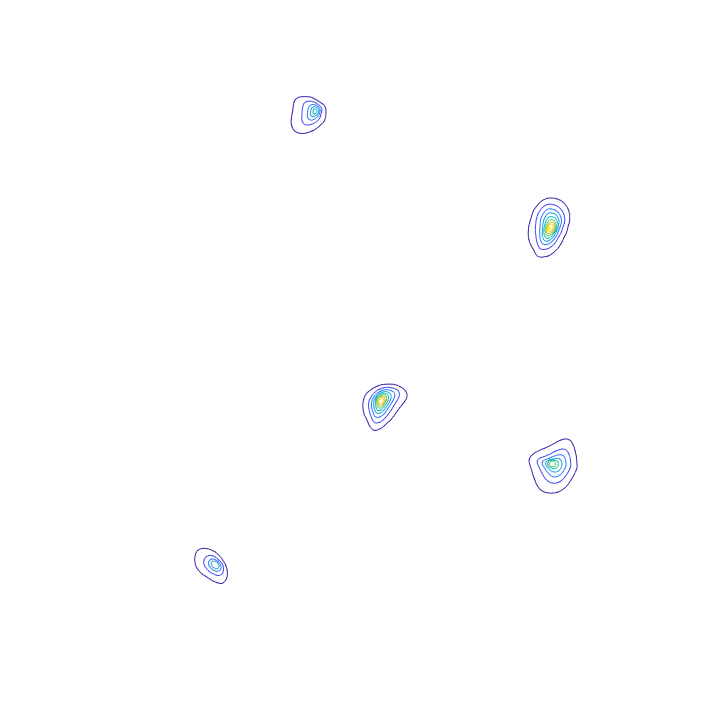}  \\
\includegraphics[height = 1.3in ]{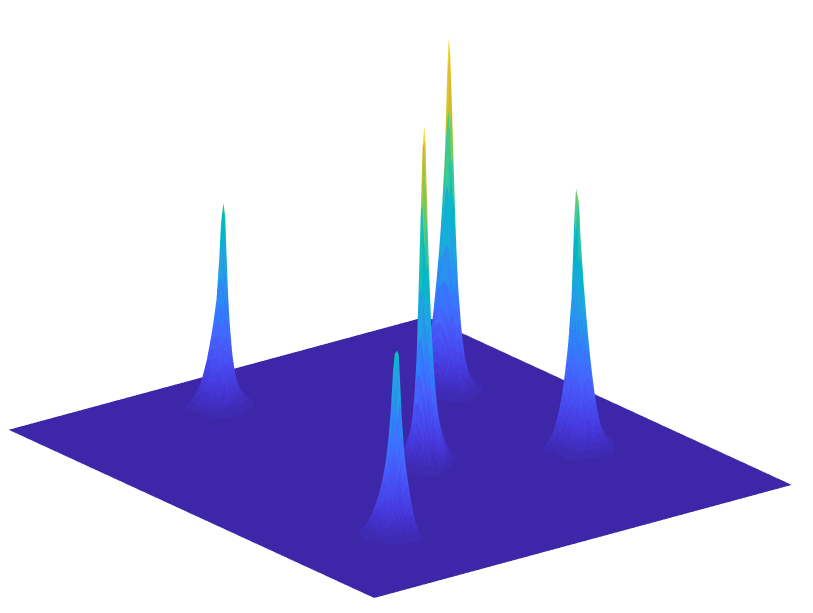}
\includegraphics[height = 1.4in ]{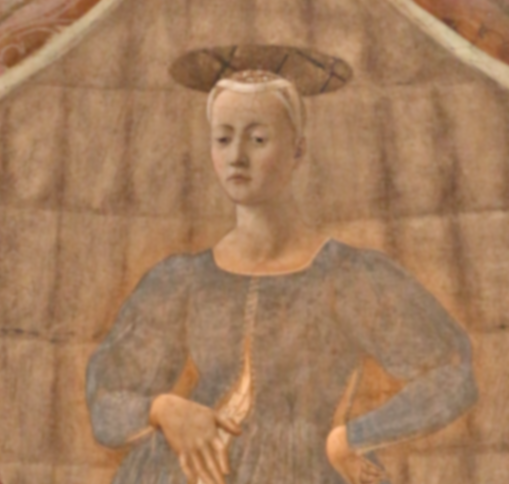} 
\caption{Heterogeneous operators in the cortex of primates: Operators are either second or fourth order, Euclidean or sub-Riemannian. They are randomly distributed with orientations prescribed by the typical orientation map of primates (pinwheels structure).
Upper left:  Ps are visualized. Upper right: Level lines of the fundamental solution. Bottom left: Surface of the fundamental solution. As in the case of rodents, the fundamental solutions are different point by point and quite symmetric. Bottom right: The reconstructed image appears free from artifacts.\newline} 
\label{f10}
\end{figure}

Ps in primates are organized in orientation maps following a pinwheels structure (Figure \ref{f10}). Also in this case there is a mixture of operators of different orders as on the previous example   but with orientations prescribed by the pinwheel structure given by Equation \eqref{pin} . Although there is a regularity in the distribution of preferred orientations, the fundamental solutions are very different from 
 the ones of homogeneous sub-Riemannian second order operators  shown in Figure \ref{f7}. In fact, the presence of a variation of orientation across the domain allows the emergence of weakly anisotropic Green functions. This feature allows a reconstruction of the perceived image without visible artifacts.\newline

In the following we present a series of simulations on different kind of images with Ps organised as in the cortex of primates. Operators are either second or fourth order, with Euclidean or sub-Riemannian metric. They are randomly distributed with orientations prescribed by the typical pinwheels structure.

\begin{figure}[htbp]
\label{f11x}
\centering
\includegraphics[height= 1.8in ]{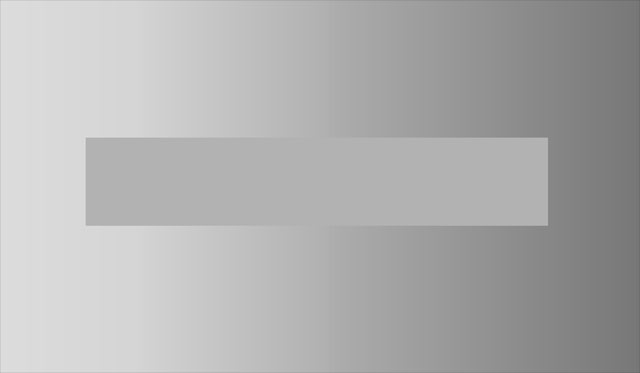}
\includegraphics[height = 1.8in ]{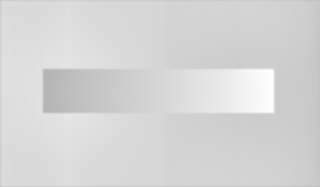}  \\
\caption{Simultaneous contrast illusion. Left: original image. The background presents a linear profile of gray and the central strip has a uniform gray value, but it is perceived as a graded gray level. Right: the outcome of simulation  shows the graded gray level typical of the perceptual illusion.\newline} 
\label{f11}
\end{figure}

In Figure \ref{f11} we considered a classical image of simultaneous contrast illusion. In the origin image the background presents a gray level increasing from left to right, and the central strip has a uniform gray value, but it is perceived as a graded gray level. The outcome of simulation  shows the perceived graded gray level in the occluding strip. 

\begin{figure}[htbp]
\label{f14x}
\centering
\includegraphics[height= 2.2in ]{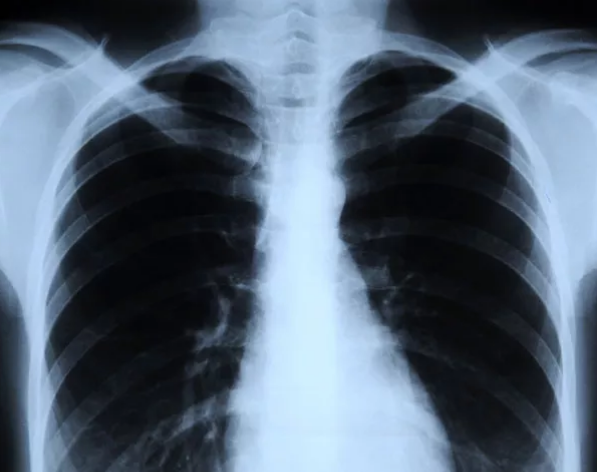}
\includegraphics[height = 2.2in ]{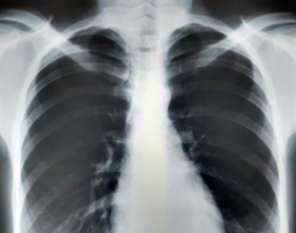}  \\
\caption{X-ray test image. Left: Original image.  Right: The cortical transform.\newline} 
 \label{f14}
\end{figure}

In Figure \ref{f14} we considered a medical image. In the reconstructed image the 
the structure of the vertebrae is much clearer. 

\begin{figure}[htbp]
\label{f12x}
\centering
\includegraphics[height= 2.9in ]{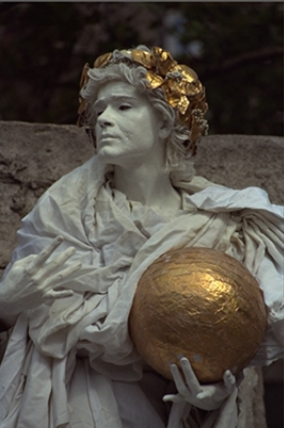}
\includegraphics[height = 2.9in ]{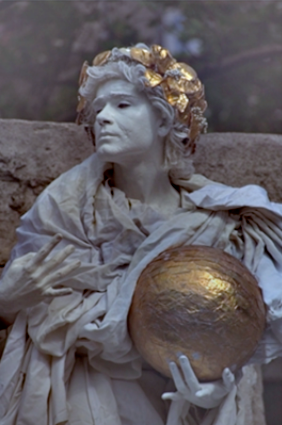}  \\
\includegraphics[height= 2.9in ]{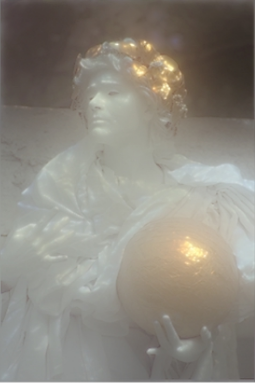}  
\includegraphics[height= 2.9in ]{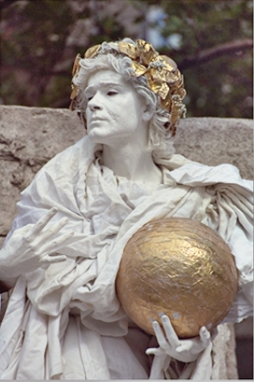}  

\caption{A test image from \cite{Ki03}. Upper Left: Original image. Upper Right: The cortical transform without any parameter. Bottom: The retinex algorithm proposed in  \cite{Ki03} with different choice of the 4 parameters.\newline} 
\label{f12} 
\end{figure}

\begin{figure}[htbp]
\centering
\includegraphics[height= 2.7in ]{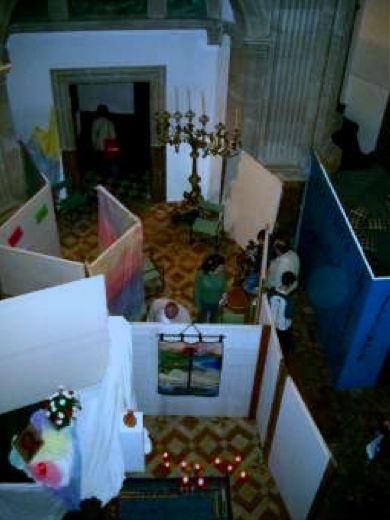}
\includegraphics[height= 2.7in ]{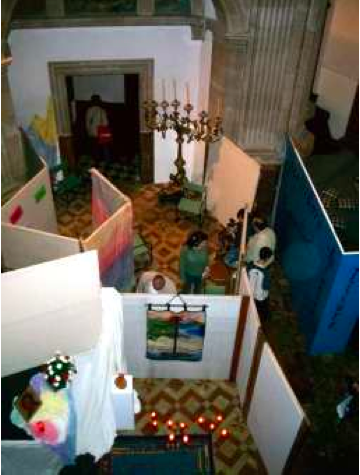}
\includegraphics[height = 2.7in ]{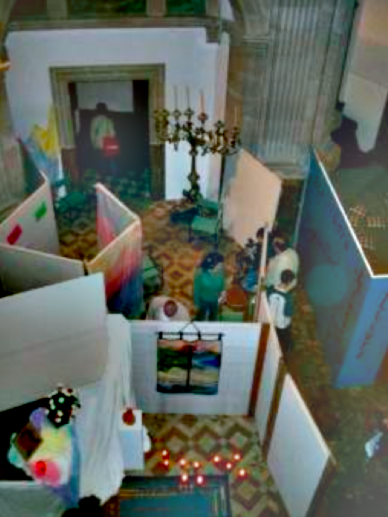}  \\
\caption{A test image from \cite{Mo11}. Left: Original image. Center: the best result of the retinex algorithm proposed in \cite{Mo11} after the choice of one parameter. Right: The cortical transform without any parameter.\newline} 
\label{f13}
\end{figure}

In Figure \ref{f12} and Figure \ref{f13} two real images are considered. The first one is taken from  \cite{Ki03}, where a retinex algorithm depending by a paramenter was proposed. Our image (top right) is comparable with the best result of their algorithm (bottom right), but the algorithm is simpler, since it does not depend on parameters. The fact that it slightly blue is compatible with perceptual laws. In facts the presence of a yellow part in the stimulus image induces  white areas to be perceived with its complementary color blue. The emergence of antagonistic pairings of colors correspond to the way our visual system is wired and it is very supported by the present model.

\section{Conclusions} 
In this paper we have presented a model of cortical transform based on the heterogeneous distribution of receptive profiles and the short range connectivity in the visual cortex. 
The action of Ps has been expressed via non negative operators with random coefficients (random metric and order), the connectivity inverts their action and the resulting cortical transform corresponds to a heterogeneous Poisson equation. 
Our main result is that the reconstruction of the perceived image is possible even in 
presence of a totally heterogeneous operator whose order varies from a point to the other. 
We provide here a formal proof of the convergence of the operator with random coefficients only in the special case of second order operator: in this case we prove convergence to a 
deterministic isotropic Laplacian. But other possibilities are at stake. For mixed order operators the experimental results are convincing even if a homogenization theory is still laking.  All these possibilities will be the object of further work in the future.

\end{document}